\documentclass[review,12pt]{elsarticle}
\usepackage{setspace}
\usepackage{subfigure}
\usepackage{multirow}
\usepackage{multicol}
\usepackage[center]{caption}
\usepackage[margin=1.3in]{geometry}
\usepackage{graphicx}
\RequirePackage{amsmath}
\usepackage{bm}
\usepackage{xfrac}
\usepackage{float}
\usepackage{arydshln} 
\usepackage{amsmath,amsfonts}

\usepackage{lineno,hyperref}
\usepackage{cleveref} 
\usepackage{nomencl}
\makenomenclature
\usepackage{etoolbox}
\renewcommand\nomgroup[1]{%
	\item[\bfseries
	\ifstrequal{#1}{A}{Acronyms}{%
		\ifstrequal{#1}{G}{Greek Symbols}{%
			\ifstrequal{#1}{O}{Other Symbols}{}}}%
	]}
\nomenclature[O]{$r$}{Euclidean distance between points (radius)}%
\nomenclature[O]{$d$}{Dimension (1D, 2D or 3D)}%
\nomenclature[O]{$s$}{Interpolated variable}%
\nomenclature[O]{$\bm{s}$}{Set of interpolated variables at the discrete points}%
\nomenclature[O]{$N$}{Total number of nodes}%
\nomenclature[O]{$\bm{A}$}{Coefficient matrix}%
\nomenclature[O]{$k$}{Degree of appended polynomial}%
\nomenclature[O]{$T$}{Temperature}%
\nomenclature[O]{$C$}{Order of convergence}%
\nomenclature[O]{$P$}{Polynomial appended to PHS-RBF}%
\nomenclature[O]{$X$,$Y$,$Z$}{Cartesian coordinates}%
\nomenclature[A]{$MQ$}{Multiquadrics}%
\nomenclature[A]{$IMQ$}{Inverse Multiquadrics}%
\nomenclature[A]{$PHS$}{Polyharmonic Splines}%
\nomenclature[A]{$TPS$}{Thin Plate Splines}%
\nomenclature[A]{$RBF$}{Radial Basis Function}%
\nomenclature[G]{$\epsilon$}{Shape parameter}%
\nomenclature[G]{$\phi$}{Radial basis function (PHS-RBF)}%
\nomenclature[G]{$\bm{\Phi}$}{Submatrix of PHS-RBF (cloud)}%
\nomenclature[O]{$\bm{P}$}{Submatrix composed of monomials appended to PHS}%
\nomenclature[G]{$\bm{\lambda}$,$\bm{\gamma}$}{Vectors containing unknown interpolation coefficients}
\nomenclature[G]{$\lambda$,$\gamma$}{Unknown interpolation coefficients}%

\journal{}

\biboptions{sort&compress}

\begin{document}
\begin{frontmatter}

\title{Application of a High Order Accurate Meshless Method to Solution of Heat Conduction in Complex Geometries }

\author{Naman Bartwal\textsuperscript{a}\fnref{Corresponding Author}}
\author{Shantanu Shahane\textsuperscript{b}}
\author{Somnath Roy\textsuperscript{a}}
\author{Surya Pratap Vanka\textsuperscript{b}}
\address{\textsuperscript{a}Department of Mechanical Engineering\\
	Indian Institute of Technology Kharagpur \\
	Kharagpur, West Bengal 721302}
\address{\textsuperscript{b}Department of Mechanical Science and Engineering\\
	University of Illinois at Urbana-Champaign \\
	Urbana, Illinois 61801}
\fntext[Corresponding Author]{Corresponding Author Email: \url{naman_bartwal@iitkgp.ac.in}}

\begin{abstract}
In recent years, a variety of meshless methods have been developed to solve partial differential equations in complex domains. Meshless methods discretize the partial differential equations over scattered points instead of grids. Radial basis functions (RBFs) have been popularly used as high accuracy interpolants of function values at scattered locations. In this paper, we apply the polyharmonic splines (PHS) as the RBF together with appended polynomial and solve the heat conduction equation in several geometries using a collocation procedure. We demonstrate the expected exponential convergence of the numerical solution as the degree of the appended polynomial is increased. The method holds promise to solve several different governing equations in thermal sciences.
\end{abstract} 
\vspace{5cm}
\end{frontmatter}

\section{Introduction}
Meshless methods have been a topic of great interest for decades due to their simplicity and further advancements.   
Some commonly used meshless techniques are smoothed particle hydrodynamics \cite{monaghan2012smoothed,ye2019smoothed,zhang2017smoothed}, reproducing kernel particle method \cite{liu1995reproducing,huang2020stabilized,patel2020meshless}, generalized finite difference method \cite{perrone1975general,liszka1980finite,gavete2017solving}, hp-clouds \cite{liszka1996hp,duarte1996hp,chen2006overview}, partition of unity \cite{melenk1996partition,babuvska1997partition,boroomand2009generalized}, element-free Galerkin \cite{belytschko1994element,abbaszadeh2020analysis,zhang2009two}, RBF-based \cite{flyer2016role,bayona2017role,bayona2019role,shahane2020high} methods, etc. Meshless methods have several advantages over commonly used grid based techniques such as unstructured finite volume and finite element methods. First and foremost is the avoidance of a grid that forms discrete elements with edges, faces, and volumes. Although the grid generation technology has significantly advanced over the years, generation of ‘good quality’ grids with desired levels of smoothness and aspect ratio is still challenging for industrial geometries. Discretization schemes are often restricted to lower order due to grid topology. Nominal second order accurate schemes can become first order in some circumstances if the grid elements are significantly skewed.  Not having a restriction on how the discrete vertices are connected gives more flexibility in interpolation of discrete values.  Second advantage with point-based methods with appended polynomial is that the local accuracy can be controlled by both the point spacing and by the appended polynomial.  This ($h-p$) refinement is more flexible than finite element and spectral element methods, where the polynomial refinement is over the entire element. It is possible to adaptively add and remove points, and also to dynamically change the local order of accuracy by recomputing coefficients only for the local cloud. Although the coefficient computations require additional work, methods to reduce this work can be developed in future.  Lastly, meshless methods provide flexibility in interface tracking, adaptive refinement, and organizing the data structures. While points still need to be generated, their flexibility in placement provides an advantage. 

Radial basis functions (RBFs) \cite{hardy1971multiquadric,fasshauer1999solving,buhmann2003radial} have been popularly used as high accuracy interpolants of function values at scattered locations over a complex domain. Common RBFs include multiquadrics (MQ), inverse multiquadrics (IMQ), Gaussians, polyharmonic splines (PHS) and thin plate splines (TPS) as given below. 
\begin{equation}
\begin{aligned}
\text{Multiquadrics (MQ): } \phi(r)=&(r^2 + \epsilon ^2)^{1/2}\\
\text{Inverse Multiquadrics (IMQ): } \phi(r)=&(r^2 + \epsilon ^2)^{-1/2}\\
\text{Gaussian: } \phi(r)=&\exp\left(\frac{-r^2}{\epsilon ^2}\right)\\
\text{Polyharmonic Splines (PHS): } \phi(r)=&r^{2a+1},\hspace{0.1cm} a \in \mathbb{N}\\
\text{Thin Plate Splines (TPS): } \phi(r)=&r^{2a} log(r),\hspace{0.1cm} a \in \mathbb{N}\\
\end{aligned}
\label{Eq:RBF_list}
\end{equation}
where, $r$ is the distance between the central point and its neighbouring data points in the domain and $\epsilon$ is the shape parameter. Several studies, including by \citet{franke1982scattered}  have evaluated the accuracy of various radial basis functions for interpolating scattered data. When appended with polynomials, the RBF interpolants can give exponential convergence as per the degree of the highest monomial appended. Since the interpolating kernels are differentiable functions, they have also been used to develop methods to solve partial differential equations governing various multi-physics phenomena. \citet{kansa1990multiquadrics} is one of the first researchers to use multiquadrics to solve partial differential equations using global interpolation over all discrete points in a complex domain. These are followed by similar works \cite{kansa1991multiquadrics,
 kansa2000circumventing,fasshauer1996solving} with applications to fluid flow, shallow water equations, electromagnetics, etc. Global methods are shown to have high accuracy but the condition numbers of the matrices to be solved are very high, making them unsuitable for large number of discrete points required to represent a complex practical geometry. Cloud-based methods, developed subsequently, interpolate the variable only around a limited number of surrounding points (clouds), leading to lower condition number of the associated matrix to be solved \cite{shu2003local}. However, the order of accuracy depends on the cloud size and the radial basis function used. Multiquadrics (MQ), inverse multiquadrics (IMQ), and Gaussians require specification of a shape parameter \cite{shu2003local,larsson2003numerical,ding2006numerical,sanyasiraju2008local,chandhini2007local,sanyasiraju2009note,vidal2016direct,zamolo2019solution} which determines the ‘flatness’ of the function and influences the accuracy as well as the condition number of the coefficient matrix. Flatter functions give higher accuracy but result in coefficient matrices of high condition numbers \cite{mairhuber1956haar,buhmann1993spectral,wong2002compactly,fornberg2004stable,fornberg2004some,larsson2005theoretical,flyer2009radial,fornberg2011stable,fasshauer2012stable,fornberg2013stable,flyer2016enhancing}. Appending polynomials to the RBF has been shown to give high accuracy, depending on the degree of the highest monomial appended \cite{barnett2015robust,flyer2016role,bayona2017role,santos2018comparing,bayona2019comparison,bayona2019role,janvcivc2019analysis,gunderman2020transport,miotti2021fully,shahane2020high,radhakrishnan2021non}.  Ideally, high order of accuracy can be achieved by appending suitably high order monomials.  However, adding more monomials also increases the size of the interpolation cloud, and the number of coefficients in the solution matrix. Fine tuning of polyharmonic splines is not required since they are not dependent on shape parameter. As shown by Fornberg and colleagues \cite{flyer2016role,bayona2017role,bayona2019role}, PHS should be appended with polynomials for stability and convergence. Moreover, this provides a powerful approach to solve partial differential equations in complex domains since the appended polynomials give higher order of accuracy. However, further research is needed to propose robust algorithms for solving problems involving industrial transport phenomena.  

Analysis of heat conduction in complex domains is important to several engineering fields including metal fabrication, materials processing, electronic cooling, additive manufacturing, cancer, retinal therapy, etc. In many situations, the concerned domains are of complex shape and only numerical solutions are feasible. In this paper, we present a systematic analysis of the PHS-RBF meshless method to solve steady and transient heat conduction in both 2D and 3D domains with emphasis on the discretization error and order of convergence. First, some simple 2D and 3D geometries are considered for error analysis. Later, three complex geometries (clamp, T-junction, and a shaft holder) are simulated to demonstrate the robustness of the PHS formulation.

\section{Methodology} \label{Sec:Method}

\begin{center}
	\begin{figure}[!htbp]
		\centering
		\includegraphics[width=0.6\textwidth]{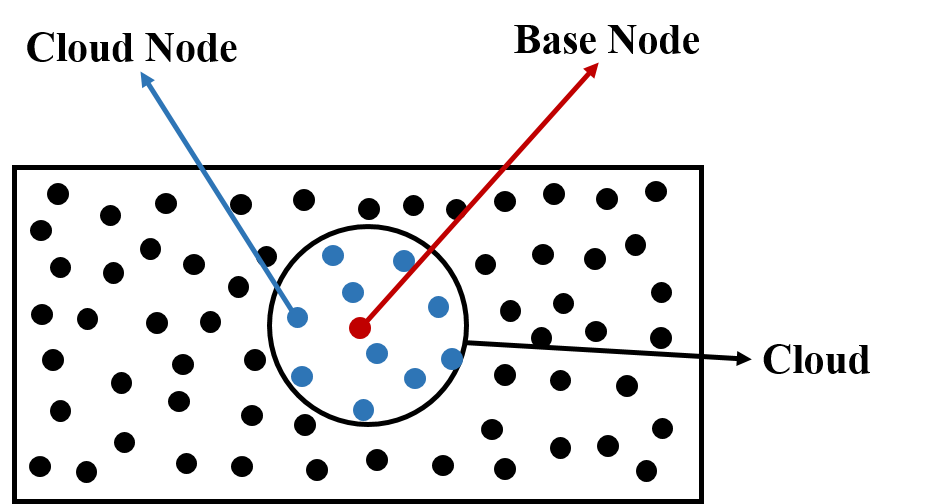}
		\caption{Nodes distribution and their associated cloud in a sample domain }
		\label{Fig:cloud_nodes}
	\end{figure}
\end{center}

In the current paper, we use a local PHS-RBF based approach with appended polynomial \cite{flyer2016role,bayona2017role,bayona2019role,shahane2020high}. \Cref{Fig:cloud_nodes} presents a sample representation of nodes distribution within a rectangular domain, where the base node (red dot) denotes the position at which the value of the function $s$ is to be evaluated in terms of the neighbouring cloud nodes (blue dots).

\vspace{-1.2cm}
 \par A continuous scalar function $s$ is approximated as a linear sum of the interpolant (polyharmonic splines with appended polynomial) over the neighbouring data points (\cref{Fig:cloud_nodes}) as shown below
\begin{equation}
s(\bm{x}) = \sum_{i=1}^{q} \lambda_{i} \phi(||\bm{x} - \bm{x_i}||_2) + \sum_{i=1}^{m} \gamma_i P_i (\bm{x})
\label{Eq:RBF_interp}
\end{equation}
\par where, $\phi \left(||\bm{x} - \bm{x_i}||_2\right)$ is a polyharmonic spline radial basis function (PHS-RBF), which is defined as $\phi(r)=r^{2a+1},\hspace{0.1cm} a \in \mathbb{N}$ and $r$ represents the Euclidean distance between the interpolation point and its neighbours in the cloud. It can be observed from \cref{Eq:RBF_interp} that there are $q+m$ unknown coefficients denoted by $\bm{\lambda}$ and $\bm{\gamma}$, which have to be determined. Additional constraints are imposed on the polynomials, resulting in $m$ equations given by
\begin{equation}
\sum_{i=1}^{q} \lambda_i P_j(\bm{x_i}) =0 \hspace{0.5cm} \text{for } 1 \leq j \leq m
\label{Eq:RBF_constraint}
\end{equation}
\par There are multiple reasons for the inclusion of these constraints such as uniquely defined solution, improvement in accuracy, etc. A detailed discussion on these constraints can be found in \cite{flyer2016role}.  

\par \Cref{Eq:RBF_interp,Eq:RBF_constraint} can be arranged in a matrix-vector system as given below.
\begin{equation}
\begin{bmatrix}
\bm{\Phi} & \bm{P}  \\
\bm{P}^T & \bm{0} \\
\end{bmatrix}
\begin{bmatrix}
\bm{\lambda}  \\
\bm{\gamma} \\
\end{bmatrix} =
\begin{bmatrix}
\bm{A}
\end{bmatrix}
\begin{bmatrix}
\bm{\lambda}  \\
\bm{\gamma} \\
\end{bmatrix} =
\begin{bmatrix}
\bm{s}  \\
\bm{0} \\
\end{bmatrix}
\label{Eq:RBF_interp_mat_vec}
\end{equation}

Here, $\bm{\Phi}$ denotes a submatrix, which is given as
\begin{equation}
\bm{\Phi} =
\begin{bmatrix}
\phi \left(||\bm{x_1} - \bm{x_1}||_2\right) & \dots  & \phi \left(||\bm{x_1} - \bm{x_q}||_2\right) \\
\vdots & \ddots & \vdots \\
\phi \left(||\bm{x_q} - \bm{x_1}||_2\right) & \dots  & \phi \left(||\bm{x_q} - \bm{x_q}||_2\right) \\
\end{bmatrix}
\label{Eq:RBF_interp_phi}
\end{equation}

 \par $\bm{\lambda}$ and $\bm{\gamma}$ represent the vectors of unknown coefficients 
 $\bm{\lambda} = [\lambda_1,...,\lambda_q]^T$, $\bm{\gamma} = [\gamma_1,...,\gamma_m]^T$, $\bm{s}$ is the set of unknowns at the discrete points i.e., $\bm{s} = [s(\bm{x_1}),...,s(\bm{x_q})]^T$. $\bm{0}$ is a submatrix of zeros of appropriate dimensions. Sizes of the submatrices $\bm{\Phi}$ and $\bm{P}$ are $q\times q$ and $q\times m$ respectively. The number of monomials appended (denoted by $m$ in \cref{Eq:RBF_interp}) for a case of given dimension ($d$) and polynomial degree ($k$) is given by $\binom{k+d}{k}$.

The unknown coefficients denoted by $\bm{\lambda}$ and $\bm{\gamma}$ can be estimated by solving \cref{Eq:RBF_interp_mat_vec}: 
\begin{equation}
\begin{bmatrix}
\bm{\lambda}  \\
\bm{\gamma} \\
\end{bmatrix} =
\begin{bmatrix}
\bm{A}
\end{bmatrix} ^{-1}
\begin{bmatrix}
\bm{s}  \\
\bm{0} \\
\end{bmatrix}
\label{Eq:RBF_interp_mat_vec_solve}
\end{equation}
\par The inverse of $\bm{A}$ in \cref{Eq:RBF_interp_mat_vec_solve} is shown only for representation. In practice, a dense direct solver is used to solve for $\bm{\lambda}$ and $\bm{\gamma}$. To perform mathematical operations such as a derivative or the Laplacian, denoted by the linear operator $\mathcal{L}$, we apply the operator $\mathcal{L}$ to the \cref{Eq:RBF_interp} as shown below.
\begin{equation}
\mathcal{L} [s(\textbf{x})] = \sum_{i=1}^{q} \lambda_i \mathcal{L} [\phi_i (\bm{x})] + \sum_{i=1}^{m} \gamma_i \mathcal{L}[P_i (\bm{x})]
\label{Eq:RBF_interp_L}
\end{equation}
\par \Cref{Eq:RBF_interp_L} can be further arranged in a matrix-vector system as follows:
\begin{equation}
\mathcal{L}[\bm{s}] =
\begin{bmatrix}
\mathcal{L}[\bm{\Phi}] & \mathcal{L}[\bm{P}]  \\
\end{bmatrix}
\begin{bmatrix}
\bm{\lambda}  \\
\bm{\gamma} \\
\end{bmatrix}
\label{Eq:RBF_interp_mat_vec_L}
\end{equation}
\par where, $\mathcal{L}[\bm{s}]$ is the differential operator evaluated at the cloud points.
\par Substituting \cref{Eq:RBF_interp_mat_vec_solve} in \cref{Eq:RBF_interp_mat_vec_L} and simplifying, we get:
\begin{equation}
\begin{aligned}
\mathcal{L}[\bm{s}] &=
\left(\begin{bmatrix}
\mathcal{L}[\bm{\Phi}] & \mathcal{L}[\bm{P}]  \\
\end{bmatrix}
\begin{bmatrix}
\bm{A}
\end{bmatrix} ^{-1}\right)
\begin{bmatrix}
\bm{s}  \\
\bm{0} \\
\end{bmatrix}
=
\begin{bmatrix}
\bm{B}
\end{bmatrix}
\begin{bmatrix}
\bm{s}  \\
\bm{0} \\
\end{bmatrix}
=
[\bm{B_1}] [\bm{s}] + [\bm{B_2}] [\bm{0}]
= [\bm{B_1}] [\bm{s}]
\end{aligned}
\label{Eq:RBF_interp_mat_vec_L_solve}
\end{equation}
\par Matrix $[\bm{B}]$ is estimated by solving $[\bm{A}][\bm{B}] = [\mathcal{L}[\bm{\Phi}] \quad \mathcal{L}[\bm{P}]]$ using a dense direct solver. Splitting $[\bm{B}]$ along columns gives the submatrices $[\bm{B_1}]$ and $[\bm{B_2}]$. $[\bm{B_1}]$ represents the coefficient matrix and multiplying it by the function values at the cloud points $(\bm{s})$ gives a numerical estimate of the differential operator evaluated at a cloud point. $[\bm{B_1}]$ depends on spatial coordinates of the cloud points and is independent of the scalar function $(s)$. Thus, $[\bm{B_1}]$ is computed and stored as a preprocessing step. These coefficients are computed for each differential operator (gradient, Laplacian etc) in the governing equations for all scattered points in the domain. The coefficients are assembled in a sparse matrix with number of rows equal to the total number of points used to discretize the entire domain. The non-zeros in every row are the coefficients of the cloud points. All the points are reordered using the RCM algorithm \cite{cuthill1969reducing} so that the neighbouring points are consecutively numbered which gives a banded sparse matrix. Governing equation discretized using the PHS-RBF coefficients is imposed on all the interior points of the domain. In this work, we focus on heat conduction problems. Hence, the governing equation has the Laplacian operator in two and three dimensions. Dirichlet boundary condition is prescribed by setting the diagonal coefficient of the row corresponding to boundary points to unity and all other entries to zero. Dot product of gradient coefficients with the outward facing normal is set as the row entries for the boundary points with Neumann condition. All the simulations in this research are performed using the open-source software MeMPhyS \cite{shahanememphys} on an i5 system with a clock frequency of 3.00 GHz.
\section{Results} 
The above algorithm is demonstrated for five 2D and 3D heat conduction problems to demonstrate accuracy and solution convergence. Apart from these, two more complex geometries are considered for transient heat conduction analyses. Among the problems, two cases have been chosen such that the numerical results can be compared directly with the exact analytical solutions. These include (i) an annular plane with the inner circle maintained at a temperature of unity, while the outer circle is at a different temperature of zero (ii) a spherical shell with similar boundary conditions. Numerical results are compared with the analytical solutions for varying degrees of appended polynomial and number of points to demonstrate the discretization accuracy of the solution procedure. Detailed investigations on the condition number and time taken to obtain the solution are also presented. Subsequently, five additional problems are studied including ellipse inside a circle (2D), spherical hole in a cuboid (3D), clamp (3D), T-junction (3D), and a shaft holder (3D). For the problems of ellipse inside a circle (2D) and spherical hole in a cuboid (3D), numerical results at different point densities and polynomial degree are compared with a benchmark solution generated using large number of points and a high degree of the polynomial.

\subsection{Steady Heat Conduction}
\subsubsection{Concentric Circles}
\label{Sec:circlecircle}
\begin{figure}[H]
	\centering
	\subfigure[Scattered points]{\includegraphics[scale = 0.49]{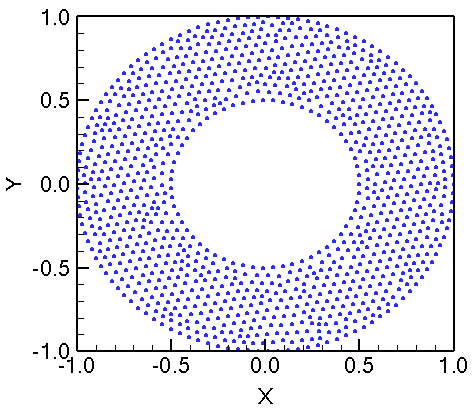}\label{fig:circlecircle_point_contour_points}}
	\subfigure[Temperature contours]{\includegraphics[scale = 0.49]{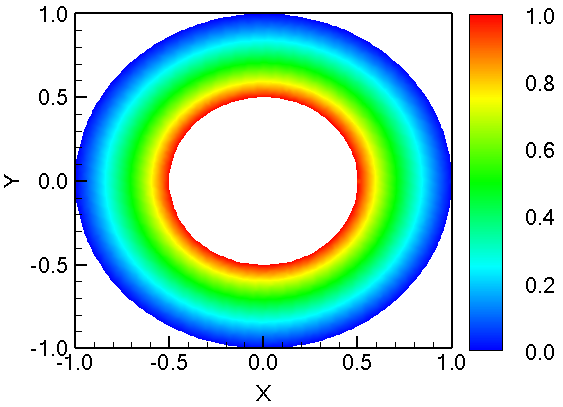}\label{fig:circlecircle_point_contour_contour}}
	\caption{(a) A sample distribution of $1073$ nodes within the domain (b) Temperature contours  }
	\label{Fig:circlecircle_point_contour}
\end{figure}
In this section, an annular geometry is considered (\cref{fig:circlecircle_point_contour_points}), where the inner circle (radius = $0.5$ unit) is maintained at a temperature of unity and the outer circle (radius = $1.0$ unit) is at a temperature of zero as shown in \cref{fig:circlecircle_point_contour_contour}. 
 Numerical solutions at varying polynomial degrees and grid sizes are compared with the exact analytical solution obtained from the conduction heat transfer equation given below.
\begin{equation}
\frac{1}{r}\frac{d}{dr}\bigg(r\frac{dT}{dr} \bigg) = 0;  
\label{Eqn:general_circle}
\end{equation} 

The temperature $T$ at any radius $r$ for the given boundary conditions is given by
\begin{equation}
T=\frac{log_{e} (r)}{log_{e} (0.5)}
\label{Eqn:ana_circle}
\end{equation}
\par
\begin{center}
 	\begin{figure}[H]  
 		\centering
 		\includegraphics[width=0.55\textwidth]{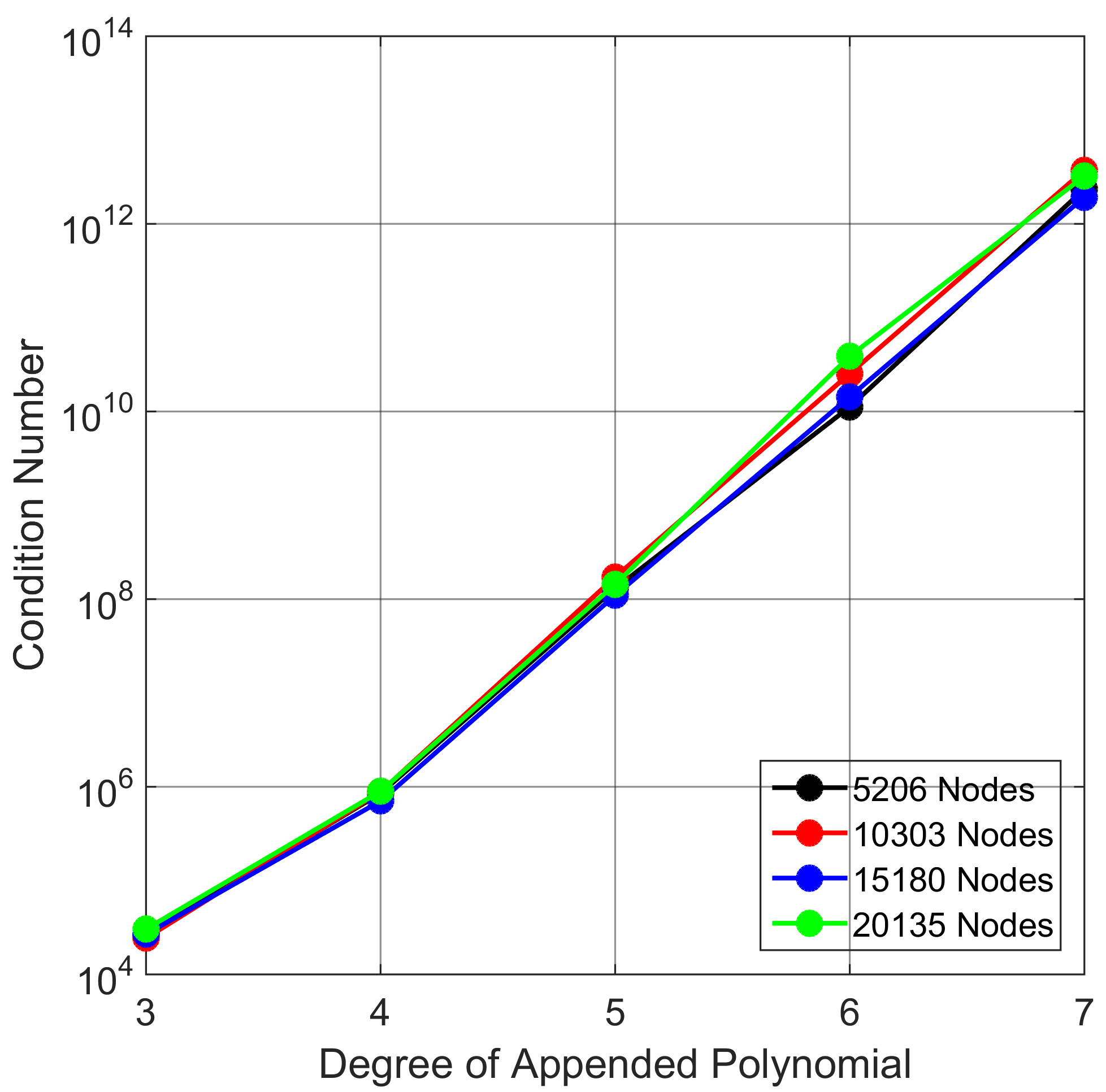}
 		\caption{Maximum condition number with the degree of appended polynomial}
 		\label{Fig:circlecircle_cond}
 	\end{figure}
 \end{center}

 Although the analytical solution is one dimensional in polar coordinates (function of $r$), we have solved in the Cartesian coordinate system, which makes the problem two dimensional. Four sets of refinements with $5206$, $10303$, $15180$, and $20135$ nodes are considered in the computational domain. \Cref{Eq:RBF_interp_mat_vec_L_solve} is solved individually for each point in the domain over its local cloud. Hence, there is a condition number of the dense linear system corresponding to each point. \Cref{Fig:circlecircle_cond} plots the maximum value of these local condition numbers over the entire domain. This is repeated for $4$ point refinements and $5$ polynomial degrees.  It can be observed that there is an exponential increase in the condition number with the degree of appended polynomial. This can be attributed to the increase in cloud size with polynomial degree, which is given by twice the number of monomials: $q = 2\binom{k+d}{k}$ \cite{shahane2020high,flyer2016role}. However, the condition number is almost independent of point refinement. This property is achieved by the technique of scaling and shifting of origin which transforms the coordinates of the local cloud to the range $[0,1]$. More details of this transformation can be found in our previous work \cite{shahane2020high}. All further calculations are carried out upto polynomial degree of six as the system becomes ill-conditioned beyond it. 
\vspace{-1.2cm}
\par
\begin{center}
 	\begin{figure}[H]  
 		\centering
 		\includegraphics[width=0.55\textwidth]{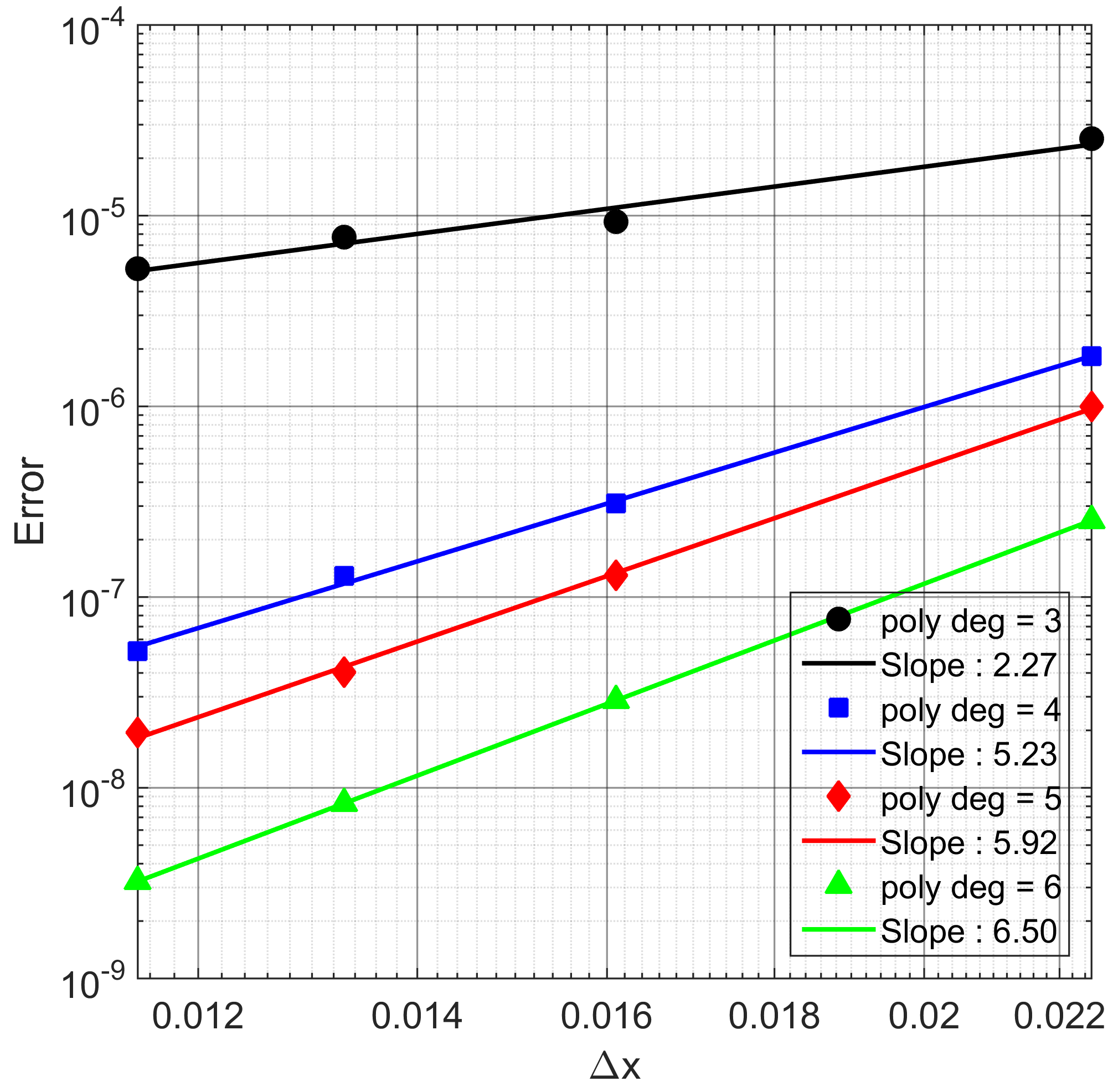}
 		\caption{Errors for concentric circles}
 		\label{Fig:circlecircle_error}
 	\end{figure}
 \end{center}

 The computed values are compared with the analytical solution of steady conduction in the annular domain (\cref{Eqn:ana_circle}). \Cref{Fig:circlecircle_error} presents the variation of the L1-norm of the error against $\Delta x$, the average distance between the nodes (computed over the entire domain). A representative value of $\Delta x$ is computed by finding the distance of the closest point for all the data points in the domain and then calculating their average.   
A line of best fit is plotted for each degree of appended polynomial. The slope of the best fit line gives an approximate estimate of the order of convergence. It can be seen from \cref{Fig:circlecircle_error} that the logarithm of the error decreases almost linearly with logarithm of $\Delta x$ for a given polynomial degree. The order of convergence is at least $k-1$ for a polynomial degree of $k$. This is consistent with the literature \cite{flyer2016role,shahane2020high}.

\subsubsection{Concentric Spheres}
\label{Sec:spheresphere}
\begin{center}
	\begin{figure}[H]  
		\centering
		\includegraphics[width=0.4\textwidth]{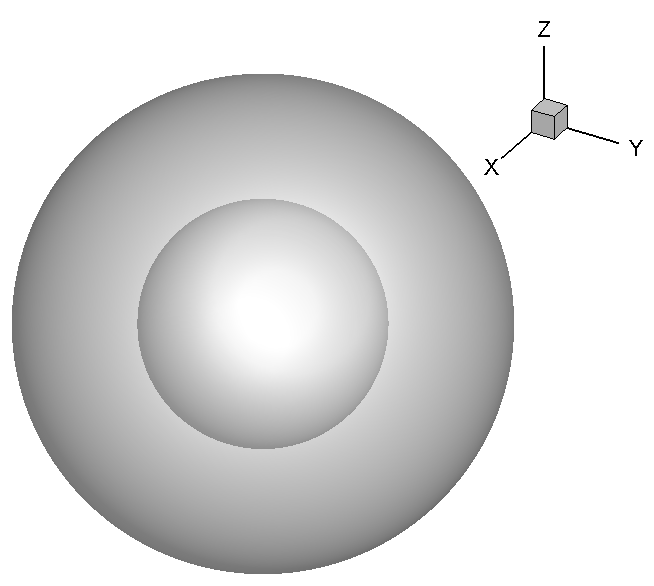}
		\caption{Geometry for concentric spheres}
		\label{Fig:spheresphere}
	\end{figure}
\end{center}
\vspace{-1.2cm}
\par To verify the solution procedure for the three dimensional case, a problem with a concentric spherical hole inside a solid sphere (\cref{Fig:spheresphere}) is considered as it can be directly compared with the exact analytical solution obtained from the heat transfer equation given below.
\begin{equation}
\frac{1}{r^2}\frac{d}{dr}\bigg( r^2\frac{dT}{dr}\bigg) = 0
\end{equation}
\par The inner sphere (radius = $0.5$ unit) is maintained at a temperature of unity and the outer sphere (radius = $1.0$ unit) is at a temperature of zero. The temperature $T$ at any radius $r$ is given by
\begin{equation}
T=\frac{1-r}{r}
\end{equation}
\vspace{-1cm}
\begin{figure}[H]
	\centering
	\subfigure[X = $0.0$ unit]{\includegraphics[width = 0.4\textwidth]{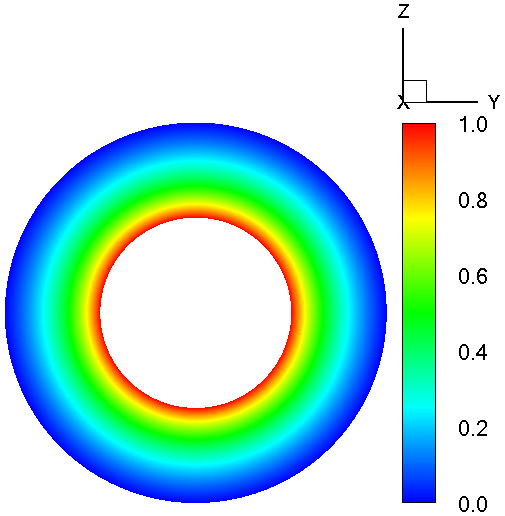}\label{fig:spheresphere_mid}}\\
	\subfigure[X = $-0.75$ unit]{\includegraphics[width = 0.34\textwidth]{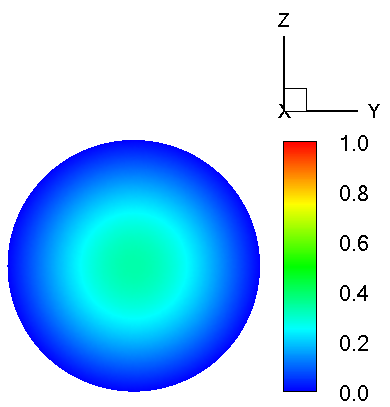}\label{fig:spheresphere_left}}
	\hspace{1cm}
	\subfigure[X = $0.75$ unit]{\includegraphics[width = 0.34\textwidth]{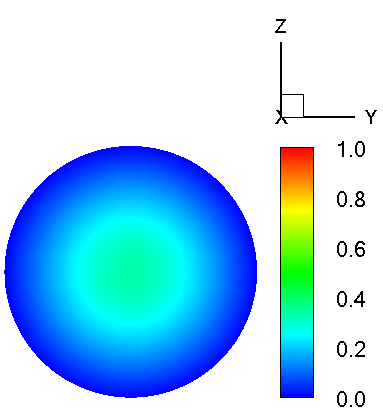}\label{fig:spheresphere_right}}
	\caption{Temperature contours at different locations parallel to YZ plane}
	\label{Fig:spheresphere_contour}
\end{figure}
\par
\begin{center}
	\begin{figure}[H]  
		\centering
		\includegraphics[width=0.57\textwidth]{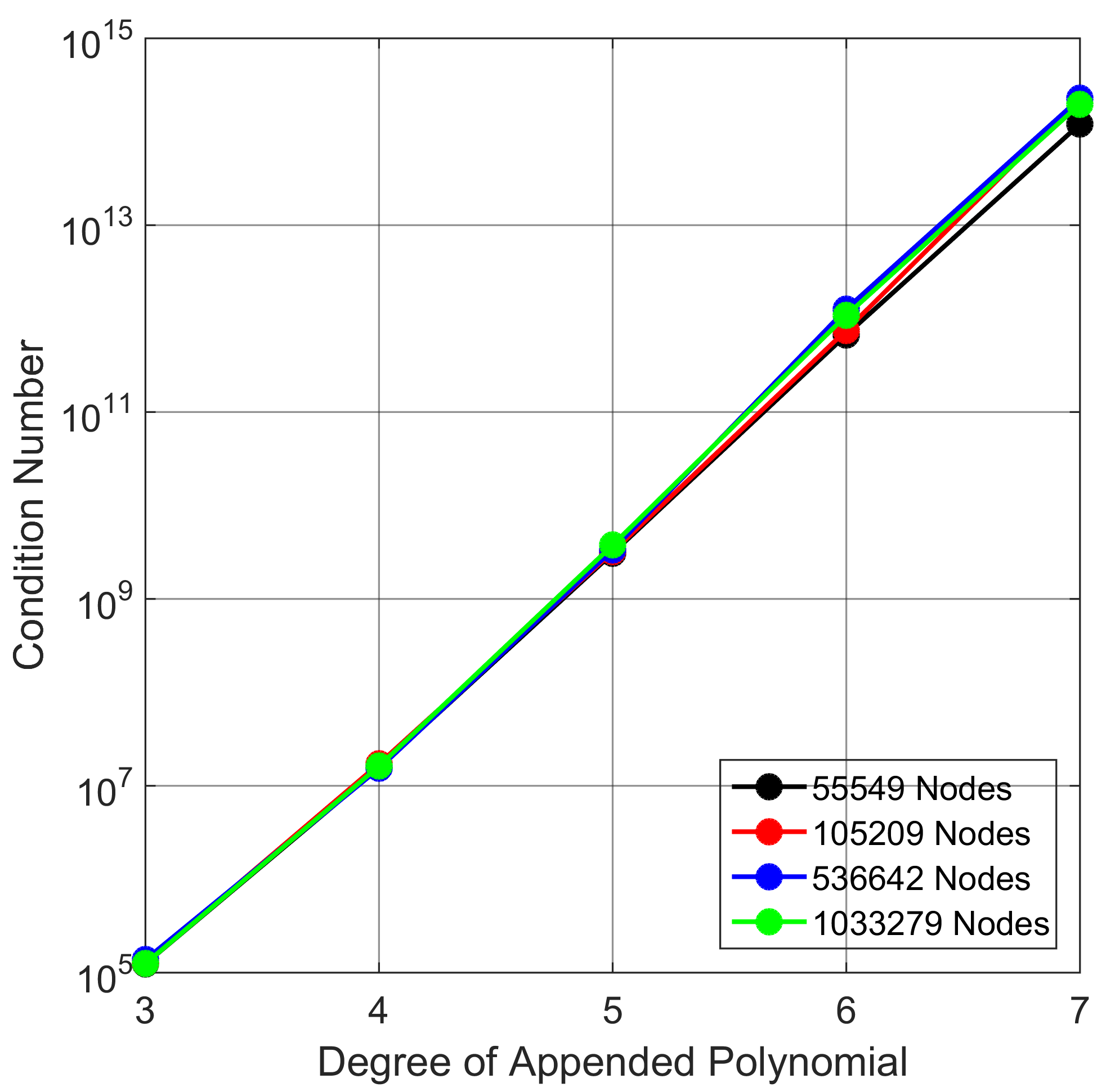}
		\caption{Maximum condition number with the degree of appended polynomial}
		\label{Fig:spheresphere_cond}
	\end{figure}
\end{center}

 Although the analytical solution is one dimensional in spherical coordinates (function of $r$), we have solved in the Cartesian coordinate system, which makes the problem three dimensional. \Cref{Fig:spheresphere_contour} depicts the temperature contours at three different cross-sections along X = $-0.75$, X = $0$, and X = $0.75$.
\Cref{Fig:spheresphere_cond} presents the condition numbers of the coefficient matrices $\bm{A}$ (\cref{Eq:RBF_interp_mat_vec_L_solve}) for this geometry with up to $1033279$ nodes. It can be seen that similar to the annular geometry case (\cref{Fig:circlecircle_cond}), the condition number varies almost exponentially with respect to the degree of the appended polynomial and can result in an ill-conditioned coefficient matrix $\bm{A}$ (\cref{Eq:RBF_interp_mat_vec_L_solve}) beyond the polynomial degree of six.
\vspace{-1.2cm}
\par
\begin{center}
	\begin{figure}[H]  
		\centering
		\includegraphics[width=0.57\textwidth]{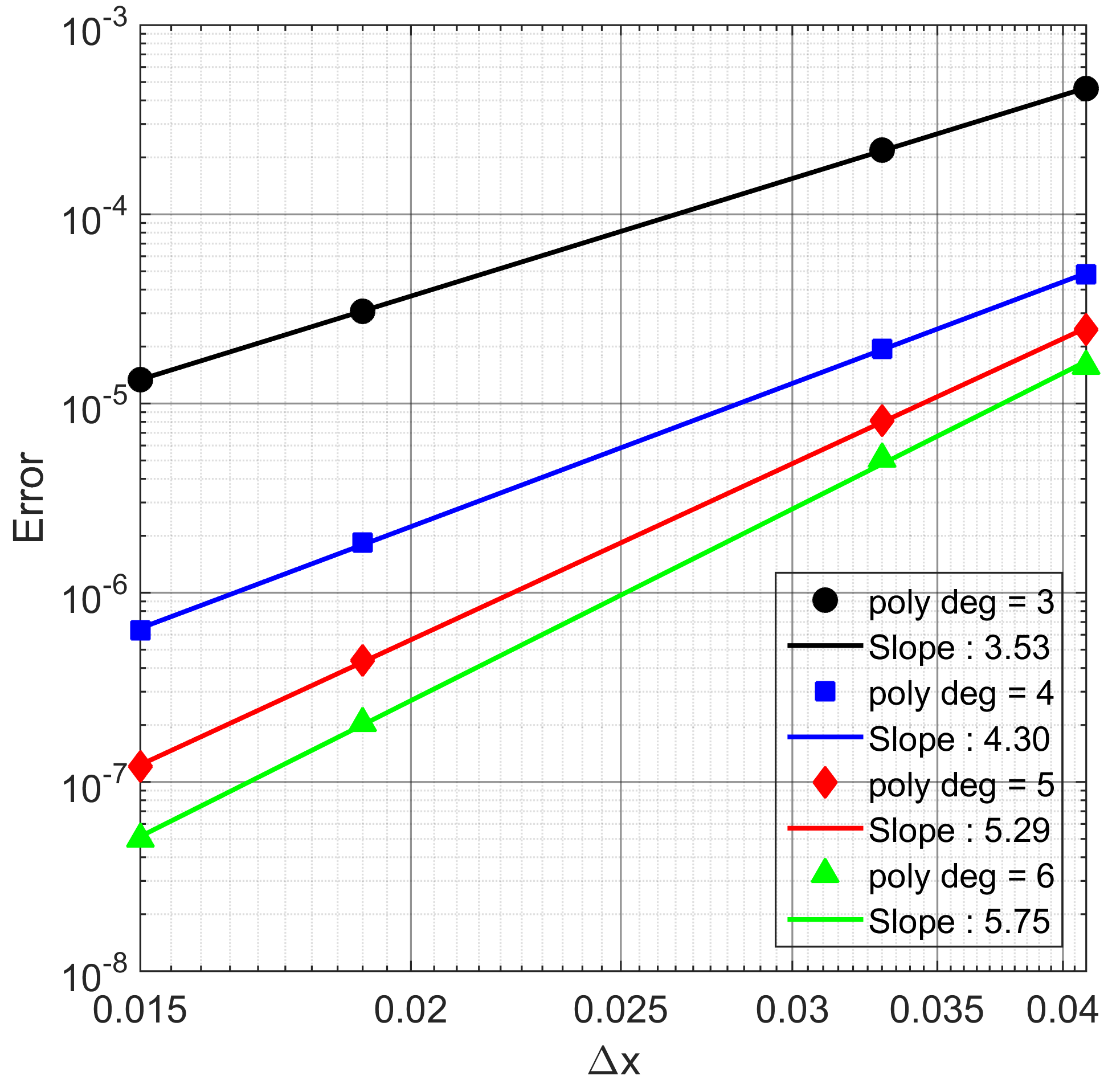}
		\caption{Errors for concentric spheres}
		\label{Fig:spheresphere_error}
	\end{figure}
\end{center}

 For comparison with the analytical solution, four sets of nodes distribution with 55549, 105209, 536642, and 1033279 nodes are considered. Similar to \cref{Fig:circlecircle_error}, it can be seen in \cref{Fig:spheresphere_error} that logarithm of the error (L1-norm) decreases almost linearly with the logarithm of $\Delta x$.  Order of convergence increases considerably with the increase in the degree of the appended polynomial and is at least $k-1$ for a polynomial degree of $k$.

\subsubsection{Ellipse inside a circle}
\begin{figure}[H]
	\centering
	\subfigure[Scattered points]{\includegraphics[scale = 0.49]{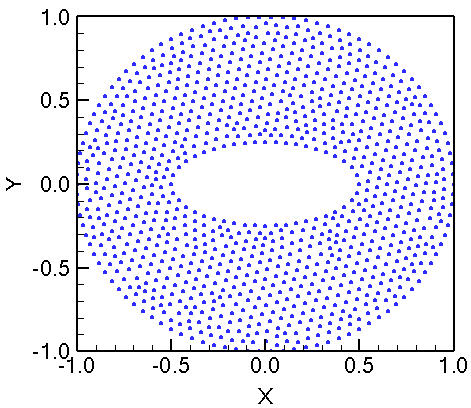}\label{fig:ellipsecircle_points}}
	\subfigure[Temperature contours]{\includegraphics[scale = 0.49]{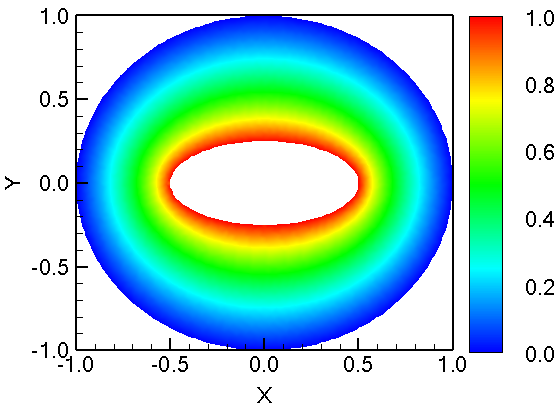}\label{fig:ellipsecircle_contour}}
	\caption{(a) A sample distribution of $1030$ nodes within the domain (b) Temperature contours  }
	\label{Fig:ellipsecircle_point_contour}
\end{figure}
\par Next, we consider a complex geometry of an ellipse embedded inside a circle. \Cref{fig:ellipsecircle_points} shows the geometry and layout of the scattered points. These points were generated from vertices of a unstructured grid generation software Gmsh \cite{geuzaine2009gmsh}. \Cref{fig:ellipsecircle_contour} shows the temperature contours for Dirichlet boundary conditions on the two perimeters (prescribed boundary conditions: inner ellipse temperature set to 1 and outer circle temperature to zero).
\vspace{1.0cm}
  \begin{table}[H]
 	\centering
 	{
 		\begin{tabular}{|l|l|l|}
 			\hline
 			\multicolumn{1}{|c|}{X} & \multicolumn{1}{c|}{Y} & \multicolumn{1}{c|}{T} \\ \hline
 			0.500 & 0.000 & \hspace{1cm}1.0 \\ \hline
 			0.583 & 0.000 & 7.24260496E-01  \\ \hline
 			0.667 & 0.000 & 5.26691444E-01  \\ \hline
 			0.750 & 0.000 & 3.66813848E-01  \\ \hline
 			0.833 & 0.000 & 2.30093193E-01  \\ \hline
 			0.917 & 0.000 & 1.09253815E-01  \\ \hline
 			1.000 & 0.000 & \hspace{1cm}0.0  \\ \hline
 			0.000 & 0.250 & \hspace{1cm}1.0 \\ \hline
 			0.000 & 0.375 & 7.66867955E-01  \\ \hline
 			0.000 & 0.500 & 5.65522307E-01  \\ \hline
 			0.000 & 0.625 & 3.92514747E-01  \\ \hline
 			0.000 & 0.750 & 2.43265416E-01  \\ \hline
 			0.000 & 0.875 & 1.13578498E-01  \\ \hline
 			0.000 & 1.000 & \hspace{1cm}0.0 \\ \hline
 		\end{tabular}
 	}
 	\caption{Reference points for ellipse inside a circle}
 	\label{tab:ellipsecircle}
 \end{table}
\newpage
\par
\begin{center}
	\begin{figure}[H]  
		\centering
		\includegraphics[width=0.57\textwidth]{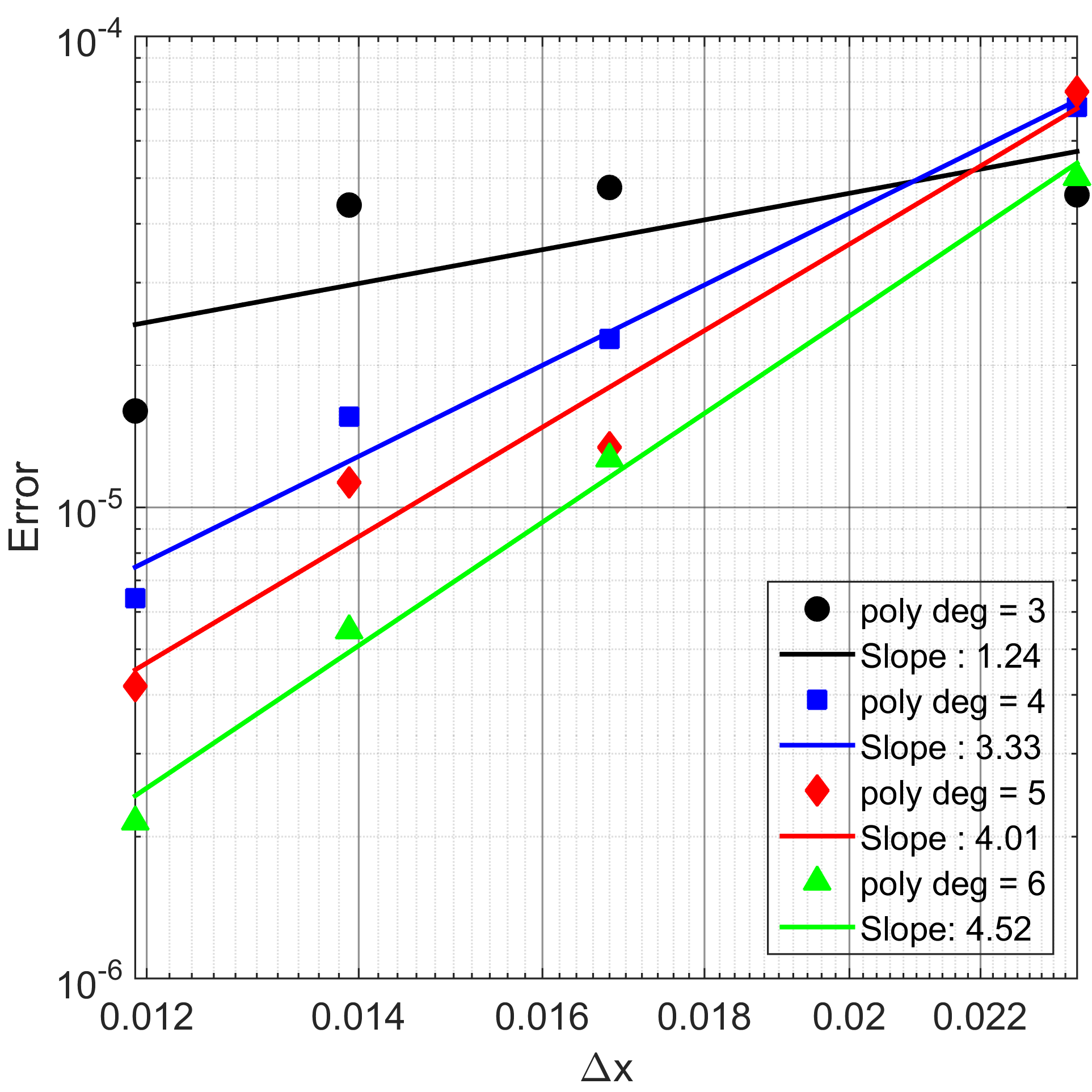}
		\caption{Errors for ellipse inside a circle}
		\label{Fig:ellipsecircle_error}
	\end{figure}
\end{center}

 To calculate the L1-norm of the discretization error, the temperatures at $14$ points along the major and minor axes of the solution domain are interpolated from a fine grid simulation (100557 points) and using a polynomial of degree six. \Cref{tab:ellipsecircle} presents the locations of the points with corresponding fine grid temperatures at these locations. The interpolated values are used as reference values and compared with results from other four sets of point distributions consisting of 5253, 10348, 15162, and 20386 points.
\Cref{Fig:ellipsecircle_error} shows the variation of the L1-norm of the error against the average $\Delta x$ on a log-log scale. As expected, a higher order of convergence is obtained as the polynomial degree is increased. Some of the slopes are slightly away from the theoretical predictions of the order of convergence. This may be due to the fact that we are estimating the error with respect to a refined numerical solution instead of an analytical solution. Moreover, the value of representative $\Delta x$ is also imprecise due to the irregular point placement.

\subsubsection{Spherical hole in a cuboid} 
\label{Sec:sphere_cuboid}
\begin{center}
	\begin{figure}[H]  
		\centering
		\includegraphics[width=0.5\textwidth]{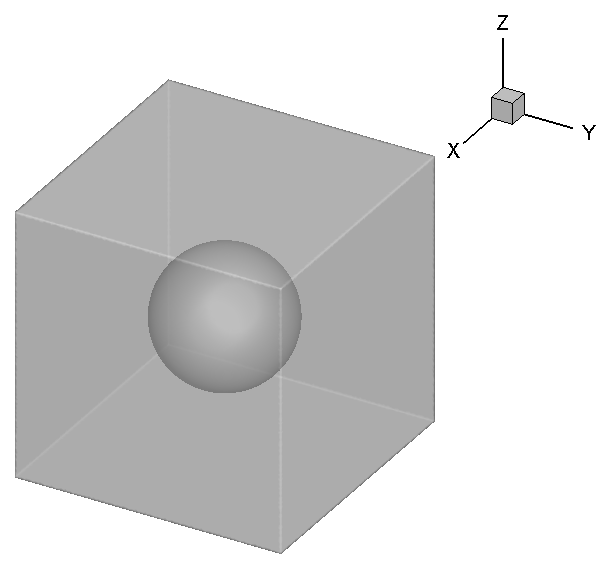}
		\caption{Geometry for spherical hole in a cuboid}
		\label{Fig:spherecuboid}
	\end{figure}
\end{center}
\vspace{-1.3cm}
\par
\begin{figure}[H]
	\centering
	\subfigure[X = $0.0$ unit]{\includegraphics[width = 0.45\textwidth ]{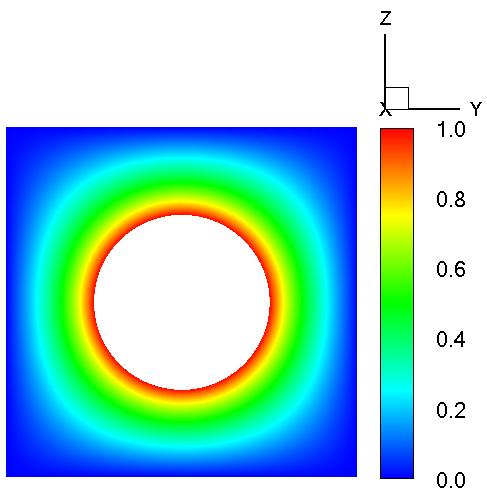}\label{fig:spherecuboid_mid}}\\
	\subfigure[X = $-0.75$ unit]{\includegraphics[width = 0.45\textwidth]{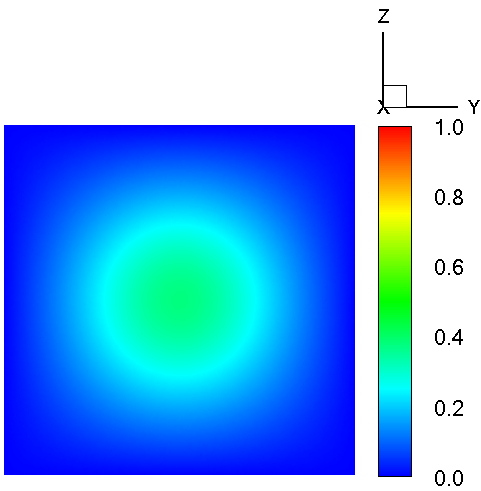}\label{fig:spherecuboid_n075}}
	\hspace{1cm}
	\subfigure[X = $0.75$ unit]{\includegraphics[width = 0.45\textwidth]{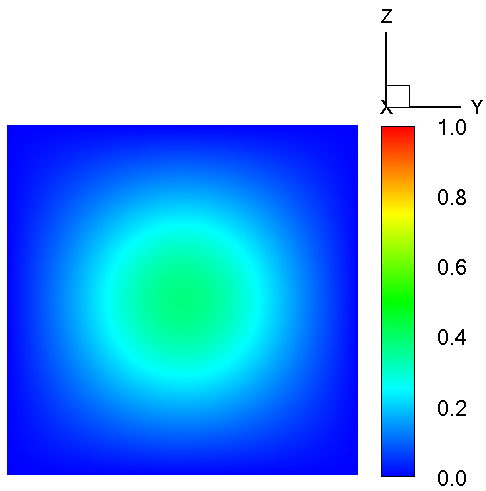}\label{fig:spherecuboid_075}}
	\caption{Temperature contours on different YZ planes}
	\label{Fig:spherecuboid_contour}
\end{figure}

 The second 3D problem considered is a spherical hole at the centre of a cuboid. The geometry is shown in \cref{Fig:spherecuboid} (cuboid : $2$ units $\times$ $2$ units $\times$ $2$ units and sphere radius = $0.5$ unit).  \Cref{Fig:spherecuboid_contour} depicts the temperature contours at three different YZ planes. The inner spherical surface is at a temperature of unity, while surfaces of the outer cuboid are maintained at a temperature of zero.  The computed fine grid solution for 1020122 nodes and a polynomial of degree six is interpolated at $14$ points along the X axis (Y = $0$, and Z = $0$). These reference temperatures (\cref{tab:spherecuboid}) are then used to evaluate the discretization accuracy by interpolating other solutions at $54873$, $103290$, and $502343$ nodes to the same discrete points.

  \begin{center}
	
	\begin{table}[H]
 		\centering
 		{
 			\begin{tabular}{|l|l|l|l|}
 				\hline 
 				\multicolumn{1}{|c|}{X}  & \multicolumn{1}{c|}{T} \\ \hline
 				0.500          & \hspace{1cm}1.0         \\ \hline
 				0.538          & 8.68951976E-01         \\ \hline
 				0.577          & 7.54963940E-01         \\ \hline
 				0.615          & 6.54505632E-01         \\ \hline
 				0.654          & 5.64896178E-01         \\ \hline
 				0.692          & 4.84045678E-01         \\ \hline
 				0.731          & 4.10285374E-01         \\ \hline
 				0.769          & 3.42252248E-01         \\ \hline
 				0.808          & 2.78807939E-01         \\ \hline
 				0.846          & 2.18980111E-01         \\ \hline
 				0.884          & 1.61918655E-01         \\ \hline
 				0.923          & 1.06861891E-01         \\ \hline
 				0.962          & 5.31095255E-02         \\ \hline
 				1.000          & \hspace{1cm}0.0        \\ \hline
 			\end{tabular}
 			
 		}
 		\caption{Reference values for spherical hole in a cuboid}
 		\label{tab:spherecuboid}
 	\end{table}
 \end{center}
\vspace{-1.2cm}
\par
\begin{center}
 	\begin{figure}[H]  
 		\centering
 		\includegraphics[width=0.57\textwidth]{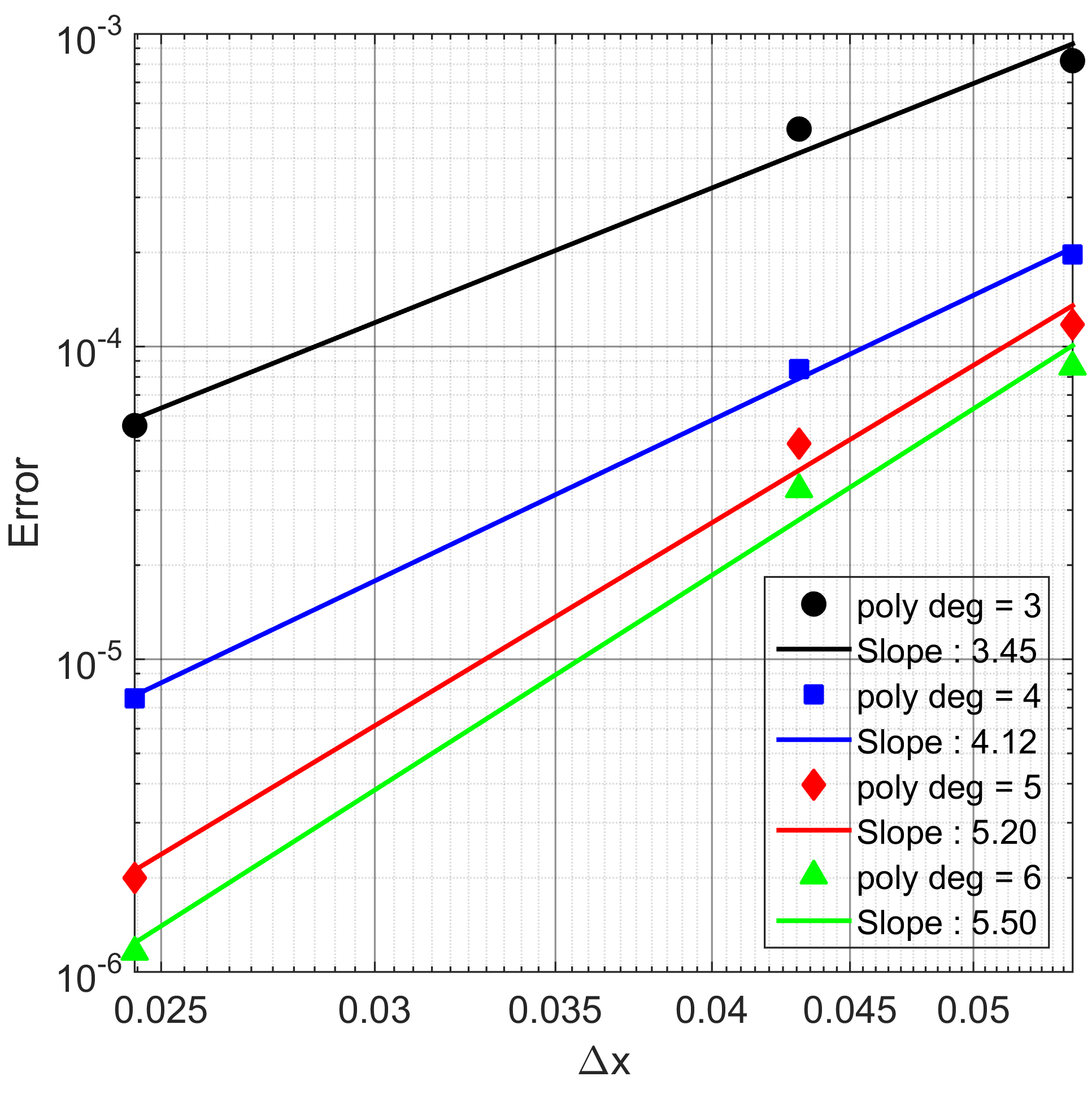}
 		\caption{Errors for spherical hole in a cuboid}
 		\label{Fig:spherecuboid_error}
 	\end{figure}
 \end{center}

 As observed in previous examples, the L1-norm of the differences decreases (\cref{Fig:spherecuboid_error}) with increase in the number of points and higher rates of convergence are obtained for higher degree of appended polynomial. The order of convergence is at least $k-1$ for a polynomial degree of $k$.
 
\subsubsection{Composite plot of observed convergence}
\par
\begin{center}
	\begin{figure}[H]  
		\centering
		\includegraphics[width=0.57\textwidth]{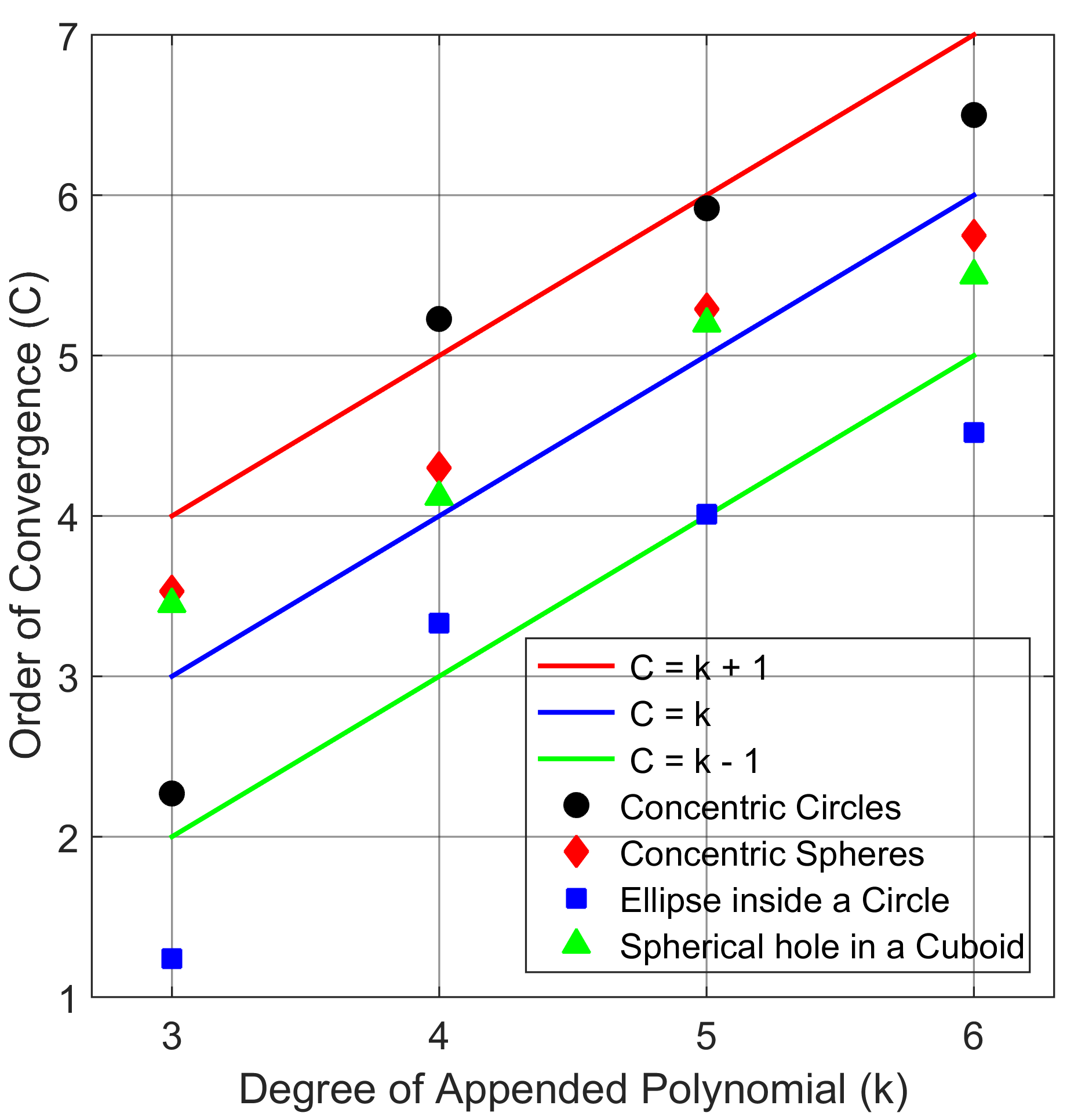}
		\caption{Order of convergence ($C$) v/s polynomial degree ($k$) for all the above discussed problems}
		\label{Fig:composite_error}
	\end{figure}
\end{center}

 In sections \ref{Sec:circlecircle} to \ref{Sec:sphere_cuboid}, we discussed the convergence of steady heat conduction equation at varying polynomial degree for both 2D and 3D problems.  As shown by \citet{shahane2020high} for a polynomial degree ($k$), the Laplacian is expected to be at least $\mathcal{O}(k-1)$ accurate. \Cref{Fig:composite_error} presents the order of convergence for different degrees of appended polynomial from $3$ to $6$. To have a better understanding of the trend of convergence, three lines ($C = k + 1$, $C = k$, $C = k - 1)$ are drawn for reference. As previously shown individually, with polynomial degree ($k$), there is a monotonic increase in the order of convergence for all the cases considered.

\subsubsection{A complex three-dimensional geometry: Clamp}
\begin{figure}[H]
	\centering
	\includegraphics[width=0.5\textwidth]{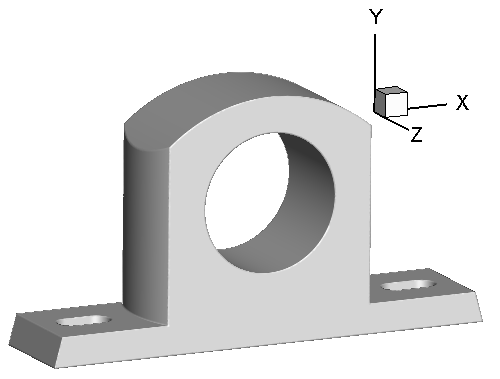}
	\caption{Geometry of a 3D Clamp}
	\label{Fig:clamp}
\end{figure}
\par
\begin{figure}[H]
 	\centering
 	\subfigure[Z = $0.0$ cm]{\includegraphics[width=0.43\textwidth]{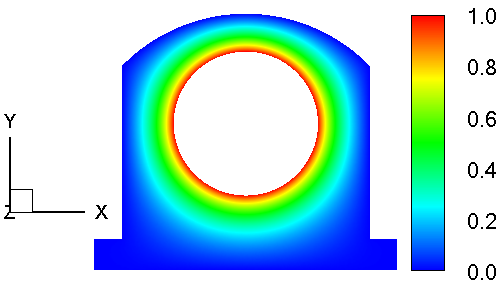}\label{fig:clamp_centre}}\\
 	\subfigure[Z = $-0.925$ cm]{\includegraphics[width=0.49\textwidth]{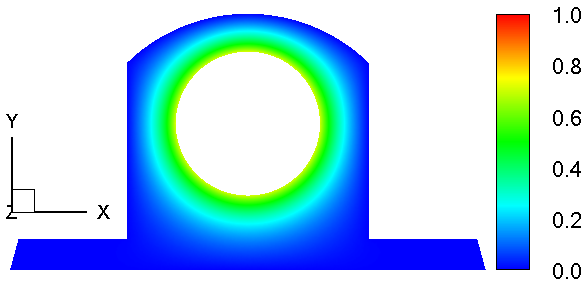}\label{fig:clamp_left}}
 	\subfigure[Z = $0.925$ cm]{\includegraphics[width=0.49\textwidth]{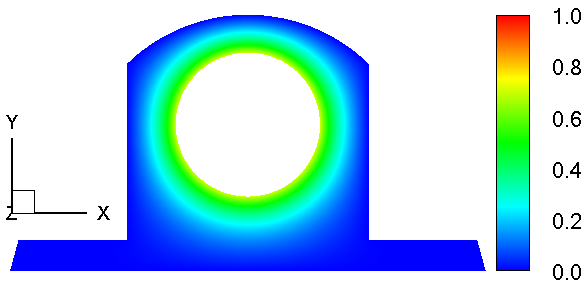}\label{fig:clamp_right}}
 	\caption{Temperature contours at different cross-sections normal to the Z-axis}
 	\label{Fig:clamp_contour}
 \end{figure}

 To investigate the robustness of the present procedure for solving heat conduction in complex shaped objects, we consider steady heat conduction in a 3D clamp of dimensions $16.5$ cm $\times$ $8.92$ cm $\times$ $3.7$ cm. A sinusoidal temperature distribution is maintained at the inner cylindrical section with temperature varying from unity at the centre to zero towards the two end faces normal to the Z-axis. All the other surfaces are maintained at a temperature of zero. Numerical simulation is carried out at $500853$ nodes and a polynomial degree of six. Temperature contours are plotted at three different cross-sections normal to the Z-axis. An expected trend of temperature distribution 
(\cref{Fig:clamp_contour}) is observed, thereby demonstrating the capability of PHS-RBF method in solving complex domains. 
\vspace{1cm}
 \newpage
\subsection{Transient Heat Conduction}
\par For checking whether the order of spatial convergence follows the same pattern in the transient case as observed in the steady heat conduction cases shown in \cref{Fig:composite_error}, we have also investigated a problem with temporal variation. Here, we consider the case of two concentric circles as presented in \cref{Sec:circlecircle} with similar Dirichlet boundary conditions using three different node refinements i.e. $9111$, $3035$, and $1073$ nodes with average $\Delta x$ of $0.0171 \, (h_{1})$, $0.0294 \, (h_{2})$, and $0.0489 \, (h_{3})$. There is almost a refinement by a factor of $1.7$. For integration in time, explicit Euler method is used, hence the temporal order of accuracy is of $\mathcal{O}(\Delta t)$. For each case, the numerically computed temperatures at a selected time of $0.4 \, s$ are interpolated at seven points along X = $0$, and seven points along Y = $0$. Richardson extrapolation is then used to obtain the spatial order of convergence ($C$) for each polynomial degree $k$. The numerically calculated value $f$ (interpolated values at three different node refinements) can be written as follows: 
\begin{equation}
\begin{aligned}
f^{h_{1}} = \alpha + \beta(h_{1})^{C} + \mathcal{O}(\Delta t) + \text{H.O.T}\\
f^{h_{2}} = \alpha + \beta(h_{2})^{C} + \mathcal{O}(\Delta t) + \text{H.O.T}\\
f^{h_{3}} = \alpha + \beta(h_{3})^{C} + \mathcal{O}(\Delta t) + \text{H.O.T}
\end{aligned}
\label{Eq:richard_interp}
\end{equation}
\par where, $h_1$, $h_2$, and $h_3$ denote the average values of $\Delta x$ and H.O.T represent higher order terms.
As the same value of $\Delta t = 1\text{E}-6$ is used for all the point distributions, we can neglect temporal errors and show that the order of convergence is given by
\begin{equation}
C = log_{1.7}\bigg(\frac{f^{h_{3}} - f^{h_{2}}}{f^{h_{2}} - f^{h_{1}}}\bigg)
\end{equation}
\par
\begin{center}
	\begin{figure}[H]
		\centering
		\includegraphics[width=0.57\textwidth]{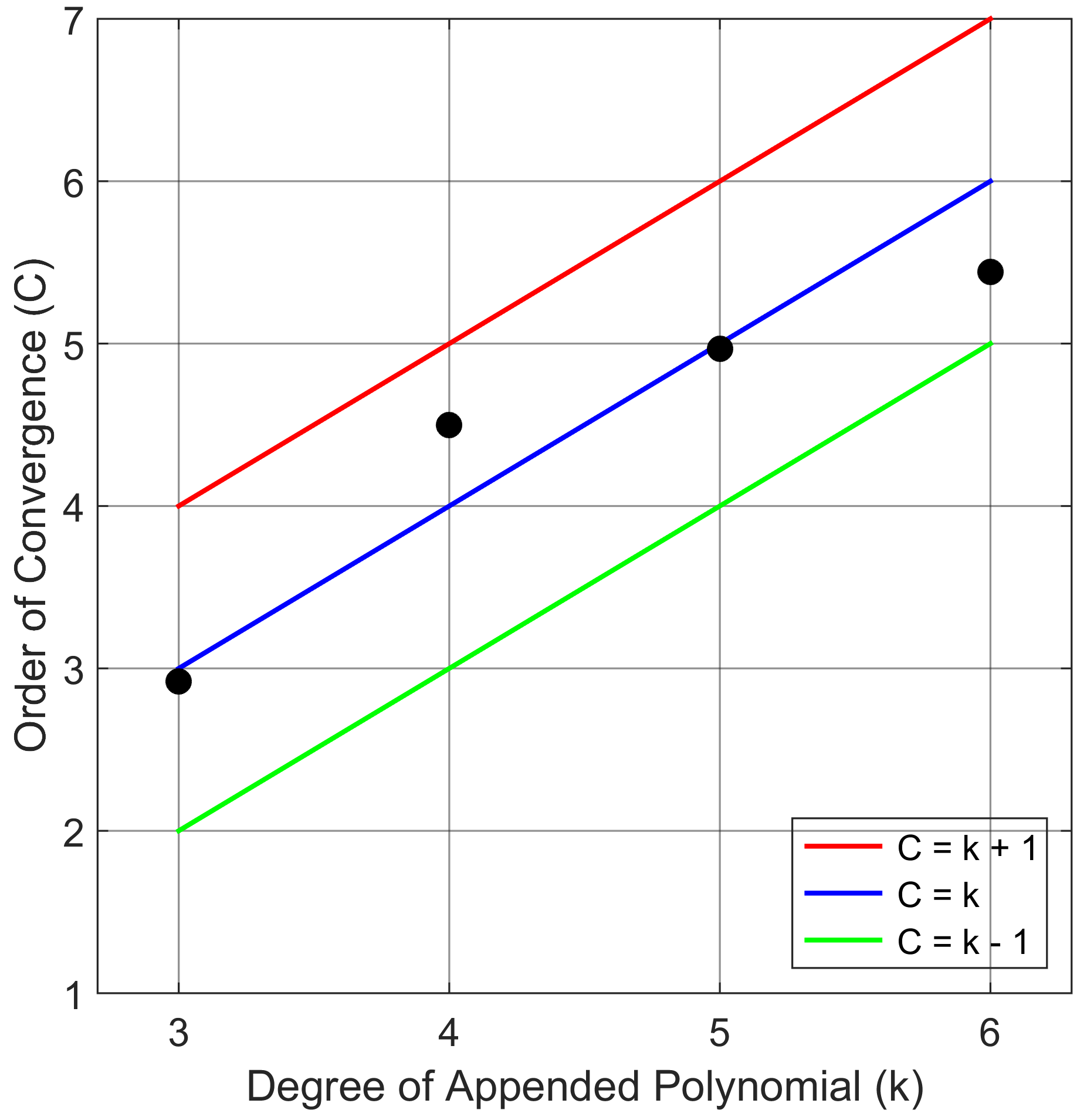}
		\caption{Order of convergence ($C$) v/s polynomial degree ($k$)}
		\label{Fig:transient_convergence}
	\end{figure}
\end{center}

 It is seen that the order of convergence increases at least as $k-1$ for a polynomial degree of $k$ (\cref{Fig:transient_convergence}), which is similar to the pattern seen in \cref{Fig:composite_error}. 
 Subsequently, two additional complex domains (T-junction and a shaft holder) are considered to demonstrate the applicability of the present meshless method to solve transient heat conduction problems in complex domains. The results, in the form of temperature contours, are presented in sections \ref{Sec:tjunc} and \ref{Sec:shaft_holder} below. 

\subsubsection{T-junction}\label{Sec:tjunc}
\begin{center}
	\begin{figure}[H]  
		\centering
		\includegraphics[width=0.5\textwidth]{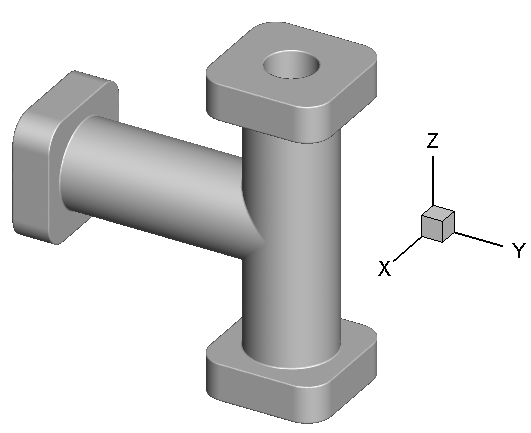}
		\caption{T-junction}
		\label{Fig:t_junc}
	\end{figure}
\end{center}
\vspace{-1.2cm}
 \par
\begin{figure}[H]
	\centering
	\subfigure[t = 0.022 s]{\includegraphics[width=0.43\textwidth]{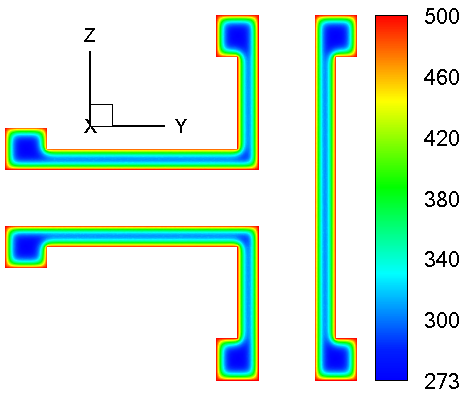}\label{fig:T_junc_first}}
	\hspace{1cm}
	\subfigure[t = 0.088 s]{\includegraphics[width=0.43\textwidth]{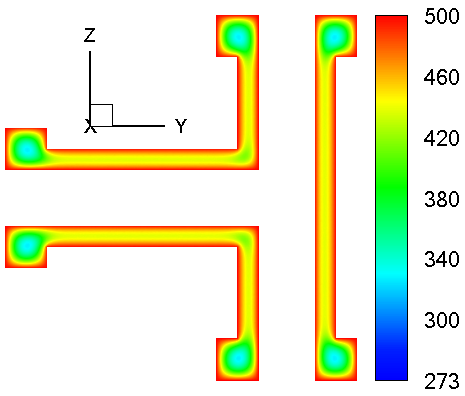}\label{fig:T_junc_second}}\\
	\subfigure[t = 0.153 s]{\includegraphics[width=0.43\textwidth]{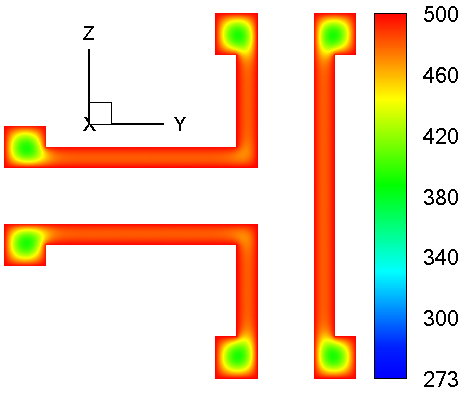}\label{fig:T_junc_third}}
	\hspace{1cm}
	\subfigure[t = 0.219 s]{\includegraphics[width=0.43\textwidth]{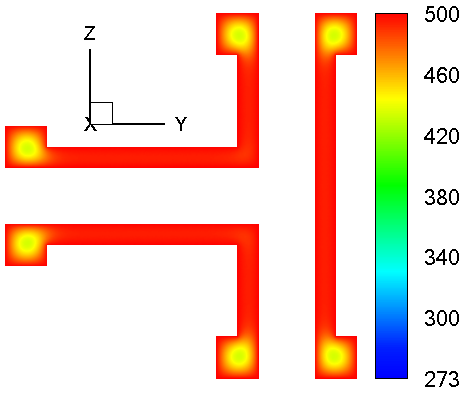}\label{fig:T_junc_fourth}}
	\subfigure[t = 0.285 s]{\includegraphics[width=0.43\textwidth]{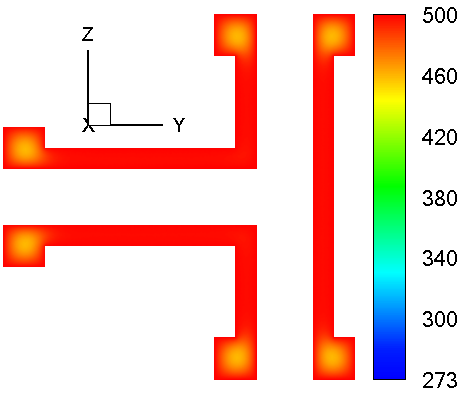}\label{fig:T_junc_fifth}}
	\hspace{1cm}
	\subfigure[t = 0.416 s]{\includegraphics[width=0.43\textwidth]{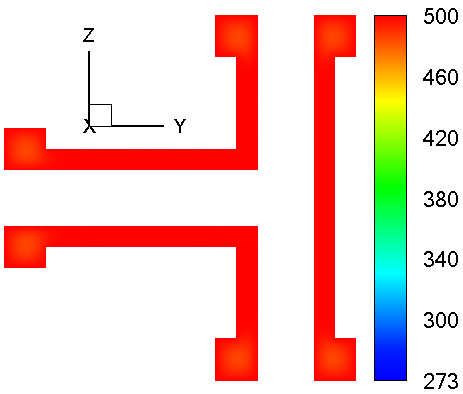}\label{fig:T_junc_sixth}}
	\caption{Temperature (K) contours at the middle cross-section parallel to the YZ plane at different time instants}
	\label{Fig:tjunc_contour}
\end{figure}

 \Cref{Fig:t_junc} shows a 3D T-junction of dimensions $5$ cm $\times$ $12.5$ cm $\times$ $13$ cm. Here, we have considered a transient heat conduction with the surfaces of the T-junction maintained at a temperature of $500$ K, while the interior is at an initial value of $273$ K.  The thermal diffusivity is taken to be $1$E$-04$ ($m^2/s$), which closely approximates the thermal diffusivity of aluminium \cite{salazar2003thermal}. The case is simulated with highly refined nodes distribution ($501422$ nodes) and a polynomial degree of six. Second order explicit Adams-Bashforth scheme is used for integration in time, and the discrete equations are accurately solved at each time step. \Cref{Fig:tjunc_contour} shows the temperature distribution on the middle cross-section parallel to the YZ plane at six different time instants.

\subsubsection{Shaft Holder}\label{Sec:shaft_holder}
\begin{center}
	\begin{figure}[H]  
		\centering
		\includegraphics[width=0.5\textwidth]{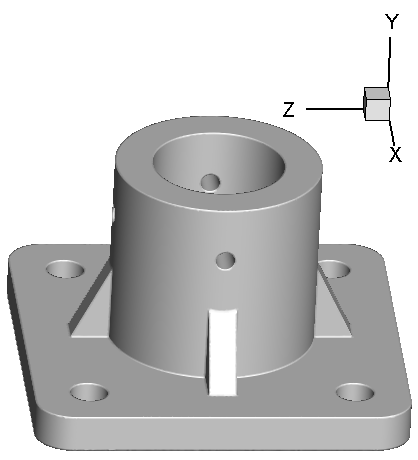}
		\caption{Shaft Holder}
		\label{Fig:shaft_holder}
	\end{figure}
\end{center}
\vspace{-1.2cm}
\par
\begin{figure}[H]
	\centering
	\subfigure[t = 0.724 s]{\includegraphics[width=0.43\textwidth]{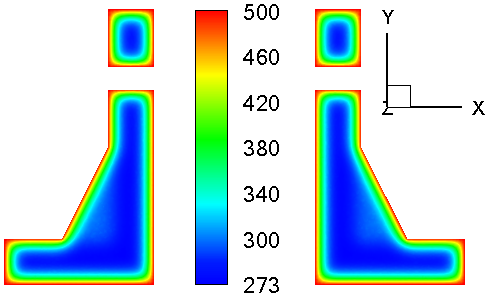}\label{fig:shaft_holder_first}}
	\hspace{1cm}
	\subfigure[t = 3.620 s]{\includegraphics[width=0.43\textwidth]{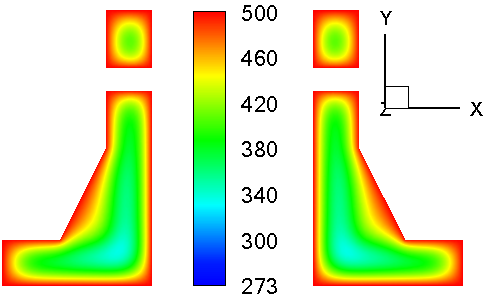}\label{fig:shaft_holder_second}}
	\subfigure[t = 6.518 s]{\includegraphics[width=0.43\textwidth]{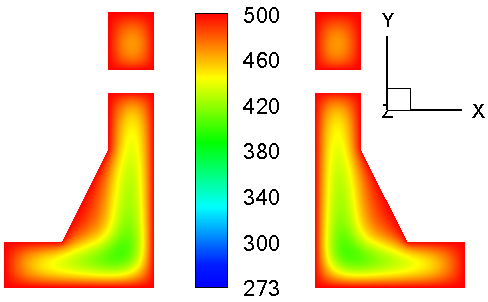}\label{fig:shaft_holder_third}}
	\hspace{1cm}
	\subfigure[t = 9.414 s]{\includegraphics[width=0.43\textwidth]{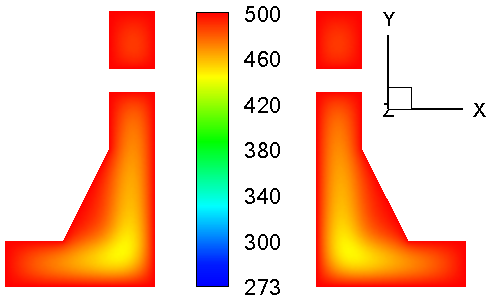}\label{fig:shaft_holder_fourth}}
	\subfigure[t = 12.311 s]{\includegraphics[width=0.43\textwidth]{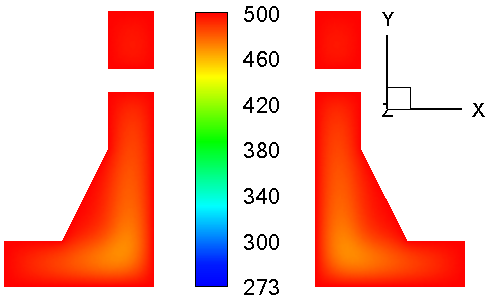}\label{fig:shaft_holder_fifth}}
	\hspace{1cm}
	\subfigure[t = 15.208 s]{\includegraphics[width=0.43\textwidth]{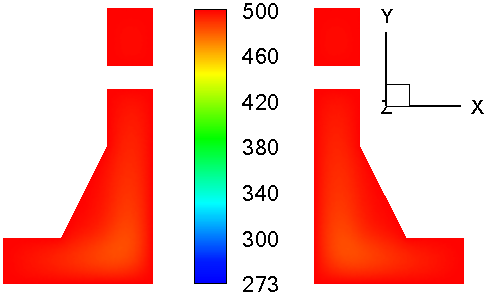}\label{fig:shaft_holder_sixth}}
	\caption{Temperature (K) contours at the middle cross-section parallel to the XY plane at different time instants}
	\label{Fig:shaft_holder_contours}
\end{figure}

 Similar to the above problem in \cref{Sec:tjunc}, another transient heat conduction problem is computed representing a 3D shaft holder (\cref{Fig:shaft_holder}) of dimensions $20$ cm $\times$ $12$ cm $\times$ $20$ cm. The temperature conditions are the same as discussed above i.e. $500$ K on all the boundaries and the initial interior value of the domain to be $273$ K. Here, we have used the thermal diffusivity of steel \cite{lienhard2005heat}, which is of the order of $1$E$-05$ ($m^2/s$). Numerical computation is carried out with $502744$ scattered points and a polynomial degree of six. Temperature contours at six different time instants are plotted on the middle cross-section parallel to XY plane as shown in \cref{Fig:shaft_holder_contours}.

\subsection{Analysis of Computational Complexity}
\par
\begin{figure}[H]
	\centering
	\subfigure[Variation with respect to the degree of appended polynomial ($k$)]{\includegraphics[width = 0.45\textwidth ]{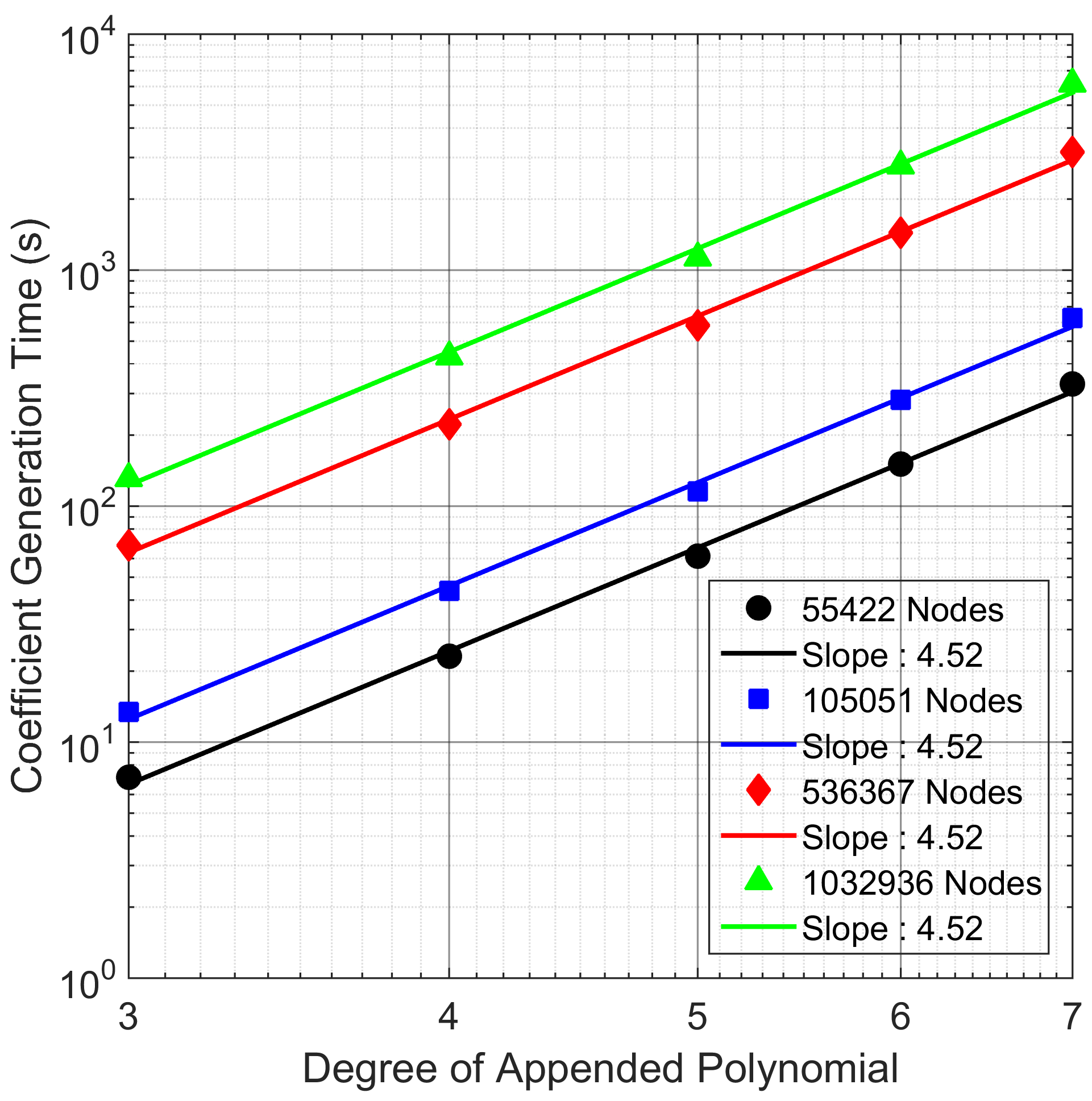}\label{fig:spheresphere_coeffgen}}
	\hspace{0.3cm}
	\subfigure[Variation with respect to the total number of nodes ($N$)]{\includegraphics[width = 0.45\textwidth ]{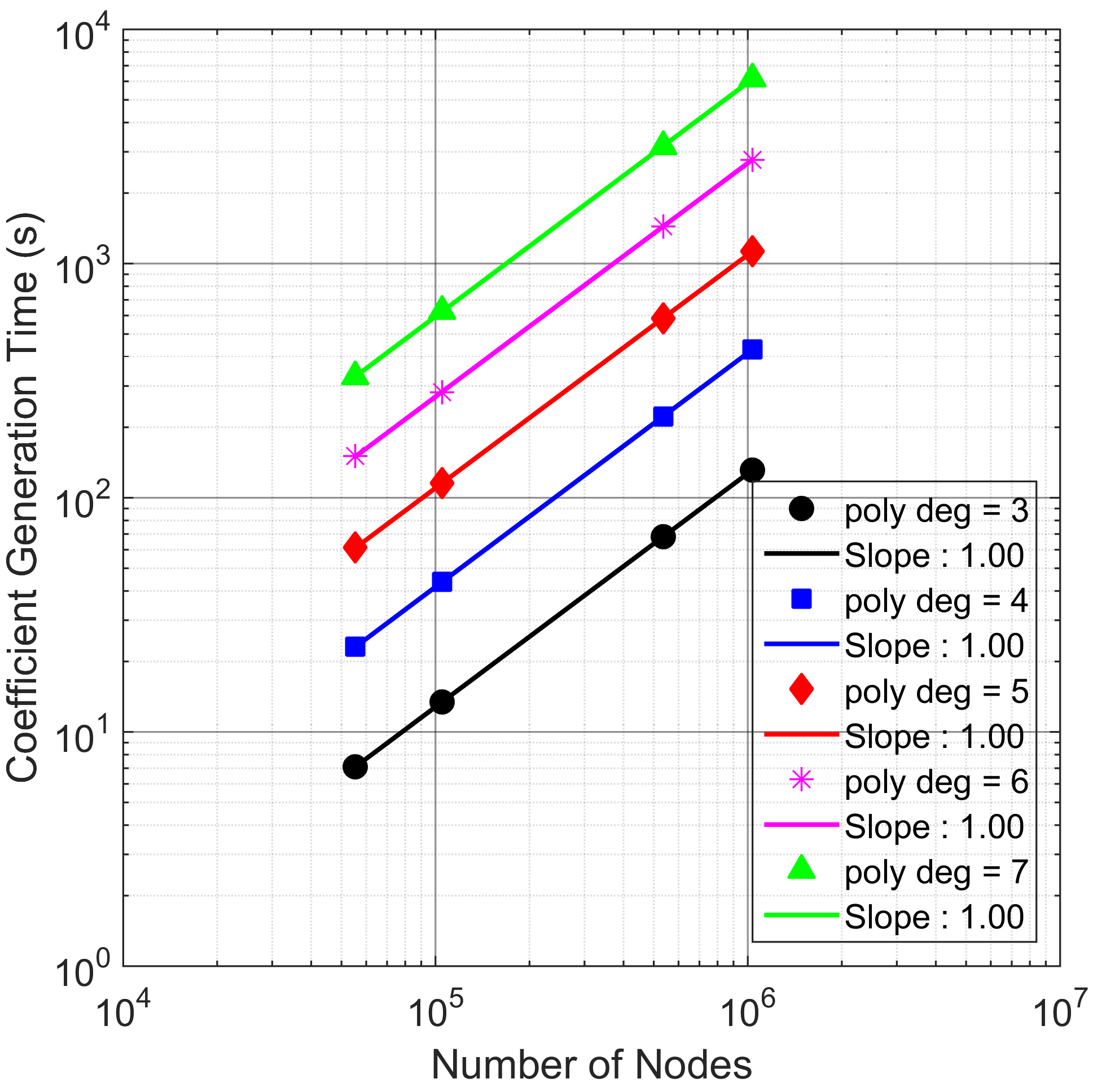}\label{fig:spheresphere_coeff_N}}
	\caption{CPU timing for coefficient generation by direct solution of a dense system}
	\label{Fig:spheresphere_coeffgen_time}
\end{figure}

\begin{figure}[H]
	\centering
	\subfigure[Variation with respect to the degree of appended polynomial ($k$)]{\includegraphics[width = 0.45\textwidth]{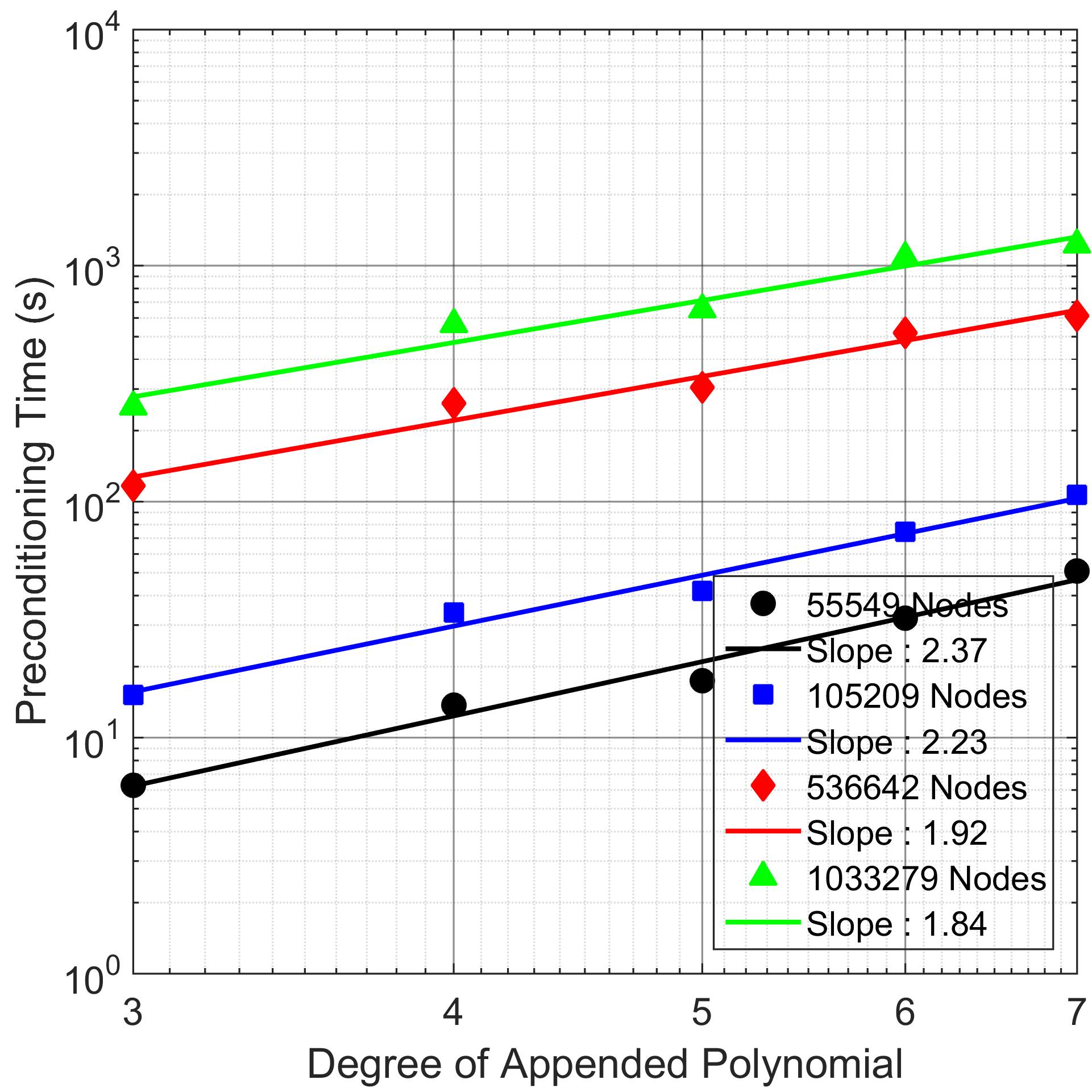}\label{fig:spheresphere_precond}}
	\hspace{0.3cm}
	\subfigure[Variation with respect to the total number of nodes ($N$)]{\includegraphics[width = 0.45\textwidth]{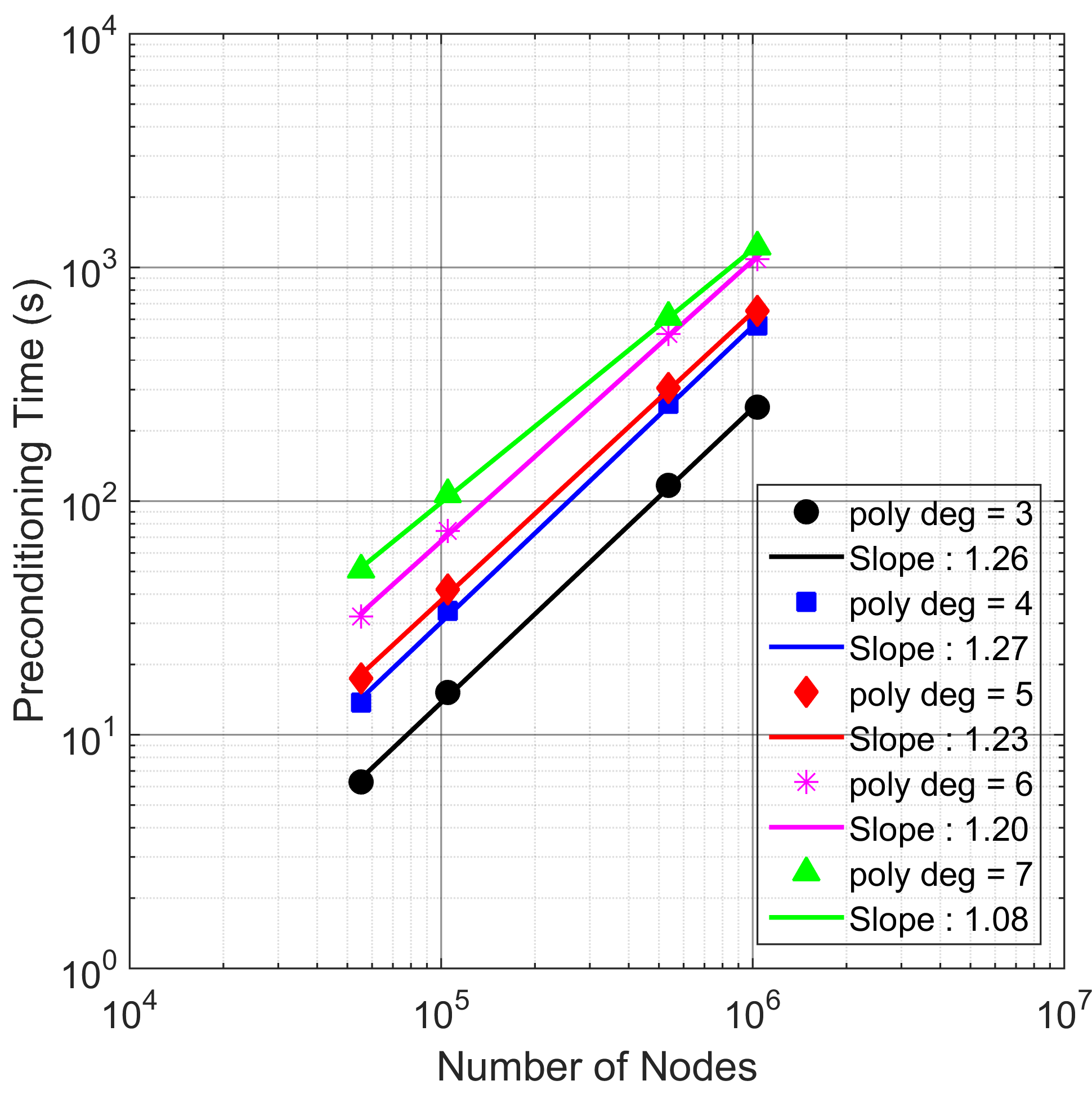}\label{fig:spheresphere_precond_N}}
	\caption{CPU timing for ILU preconditioning}
	\label{Fig:spheresphere_precond_time}
\end{figure}

\begin{figure}[H]
	\centering
	\subfigure[Variation with respect to the degree of appended polynomial ($k$)]{\includegraphics[width = 0.45\textwidth]{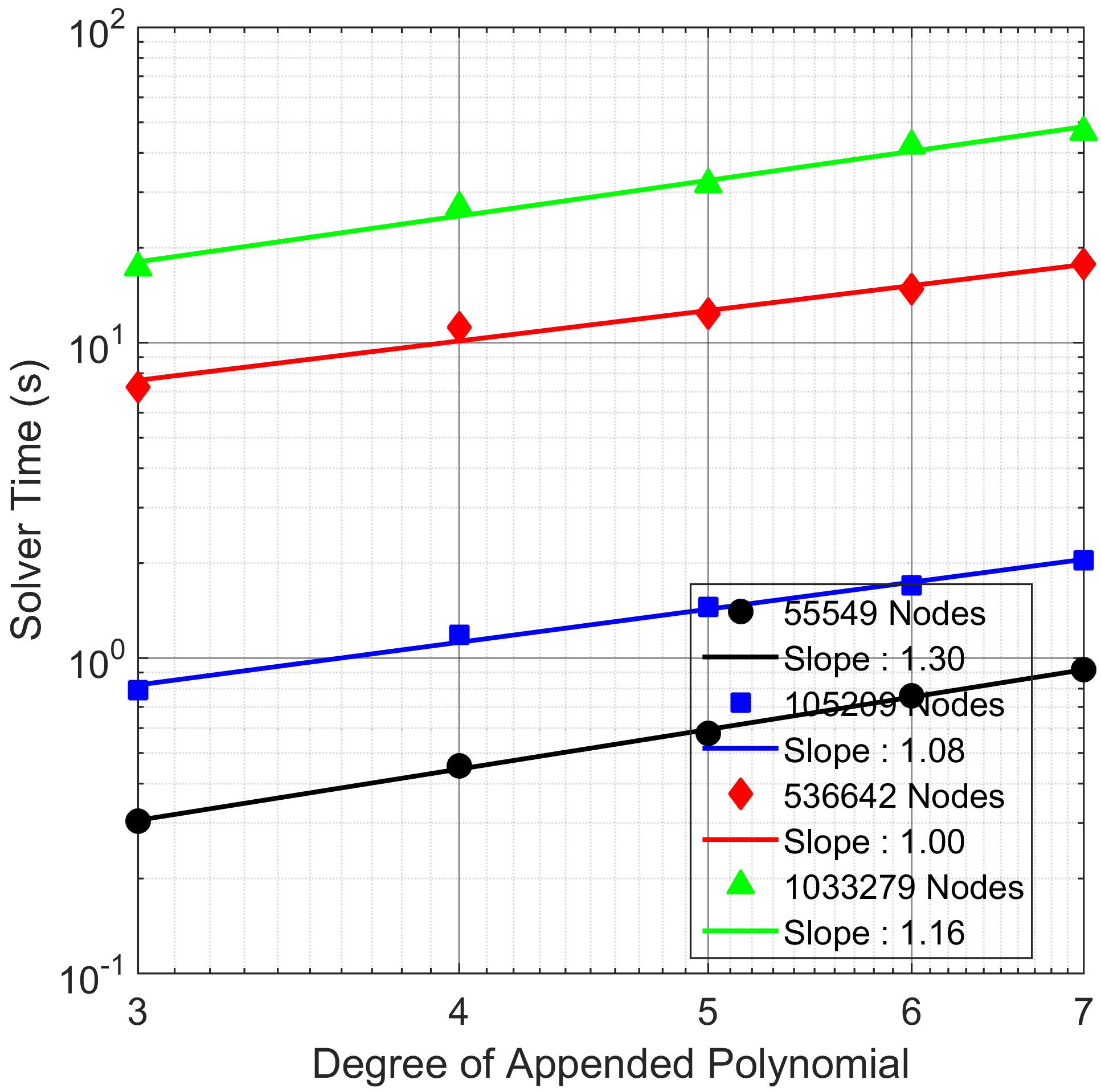}\label{fig:spheresphere_solver}}
	\hspace{0.3cm}
	\subfigure[Variation with respect to the total number of nodes ($N$)]{\includegraphics[width = 0.45\textwidth]{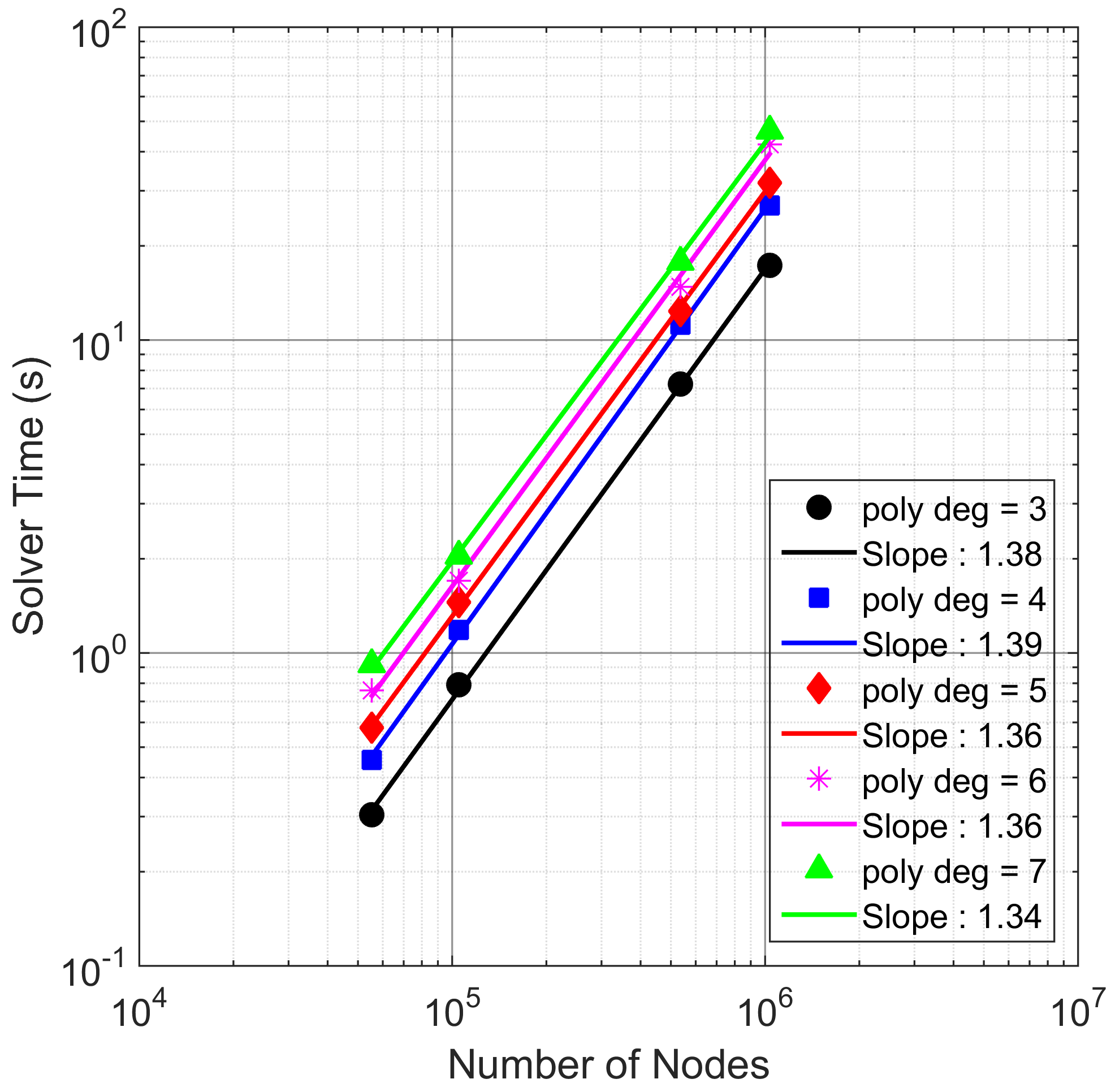}\label{fig:spheresphere_solver_N}}
	\caption{CPU timing for sparse BiCGSTAB solver}
	\label{Fig:spheresphere_solver_time}
\end{figure}

 Figures \ref{Fig:spheresphere_coeffgen_time}\textendash\ref{Fig:spheresphere_solver_time} present times for different steps in the computation. These timings correspond to heat conduction in concentric spheres (\cref{Sec:spheresphere}), but similar trends will be seen for other problems. \Cref{Fig:spheresphere_coeffgen_time} presents the coefficient generation times corresponding to the direct solution of the dense system in \cref{Eq:RBF_interp_mat_vec_L_solve}. For each point in the domain, the size of the dense matrix $\bm{A}$ for a polynomial degree of $k$ in three dimensional problems is given by $3\binom{k+d}{k}  = \frac{(k+1)(k+2)(k+3)}{2}$. \Cref{fig:spheresphere_coeffgen} shows that the computational cost is proportional to $k^{4.52}$. Since this operation is performed for each point in the domain, the cost is linear in the total number of points $N$ (\cref{fig:spheresphere_coeff_N}). Thus, we see that the computational cost is of $\mathcal{O}(Nk^{4.52})$. However, this CPU time can be reduced by grouping the points into clusters and sharing the cloud between the points in the cluster. Further the points can be processed in parallel on multiple cores, since the operations at any point are independent from the other points.
\par \Cref{Fig:spheresphere_precond_time} shows the preconditioning times. Here we have used an ILU preconditioning with BiCGSTAB as the solver. The bandwidth of the final assembled sparse coefficient matrix is equal to the cloud size given by $2\binom{k+d}{k} = \frac{(k+1)(k+2)(k+3)}{3}$. Hence, the total number of non-zeros is of $\mathcal{O}(Nk^{3})$ for 3D problems. From figs. \ref{fig:spheresphere_precond} and \ref{fig:spheresphere_precond_N}, we see that the computational cost of ILU preconditioning is approximately quadratic in $k$ and slightly superlinear in $N$. The reason for the superlinear increase with $N$ may be that for a fixed $k$, number of non-zeros in the preconditioning matrix may be increasing with $N$. Essentially, large matrices with same initial bandwidth need a stronger preconditioner with higher bandwidth. However, since both coefficient generation and preconditioning are pre-processing steps, their CPU time can be amortized over the time marching steps in case of transient problems.
\par \Cref{Fig:spheresphere_solver_time} shows the solver time. Here BiCGSTAB is used to solve the ILU-preconditioned sparse matrix.  From figs. \ref{fig:spheresphere_solver} and \ref{fig:spheresphere_solver_N}, we observe that the CPU time is slightly superlinear with respect to both $k$ and $N$.
 \par We have currently used the C++ template library Eigen for solutions of both the dense and sparse systems as well as preconditioning \cite{eigenweb}. Parallelization over multiple CPU or GPU cores can reduce these computational times.

\section{Summary and Future Work}
\par The current work applies a high-accuracy radial basis function (RBF) meshless method to compute heat conduction in complex domains. The method uses scattered distribution of points to interpolate discrete temperatures using polyharmonic splines with appended polynomial. First, two problems have been chosen so that the steady state numerical results can be directly compared with analytical solutions, validating the numerical technique and assessing the accuracy. Subsequently, more complex problems are considered for which benchmark solutions are obtained with large numbers of nodes and high degree of appended polynomial. These results are then used for comparing with results obtained with less numbers of points and varying degree of appended polynomial. These problems demonstrate the applicability and accuracy of the PHS-RBF to heat transfer in complex domains. All the computations are carried out upto a polynomial degree of six. It is seen that for coefficient generation, the computational cost is proportional to $k^{4.52}$, approximately quadratic in $k$ for ILU preconditioning, and slightly superlinear with respect to $k$ for solving the ILU-preconditioned sparse matrix. It is shown that for a given number of points, the discretization error (L1-norm) decreases rapidly as the degree of appended polynomial is increased. The order of convergence is found to be at least $k-1$ for a polynomial degree of $k$ for both steady and transient cases. An expected trend of temperature distributions is observed in two transient problems, thus demonstrating the capability of the method to solve transient heat conduction problems in complex domains. Future extensions of the method include adaptive local refinements in point spacing as well as degree of appended polynomial. This will permit controlling accuracy in regions of steep variation.

\printnomenclature

\section{Declaration of Originality}
The authors declare that this is a novel and original work and it has not been published elsewhere and is not currently being considered for publication elsewhere. All tables and figures are work of the authors and no permissions are required.

\bibliography{References}

\end{document}